\documentclass[11pt]{article} 
\usepackage{amsmath,amssymb,amsfonts,amsthm,amscd,graphicx,psfrag,epsfig,mathrsfs}
\usepackage{color}
\definecolor{blue}{rgb}{0,0,0.7}
\definecolor{red}{rgb}{0.75, 0, 0}
\usepackage{stmaryrd}
\usepackage[titletoc,toc]{appendix}
\usepackage{hyperref}
\hypersetup{colorlinks, linkcolor=blue}
\usepackage[all]{xy}           
\usepackage{marginnote}  
\usepackage{float}

\newtheorem{theorem}{Theorem}[section]

\newtheorem{theorem-definition}[theorem]{Theorem-Definition}
\newtheorem{theorem-construction}[theorem]{Theorem-Construction}
\newtheorem{lemma-definition}[theorem]{Lemma--Definition}
\newtheorem{proposition-definition}[theorem]{Proposition--Definition}
\newtheorem{lemma}[theorem]{Lemma}
\newtheorem{proposition}[theorem]{Proposition}
\newtheorem{corollary}[theorem]{Corollary}
\newtheorem{conjecture}[theorem]{Conjecture}
\newtheorem{definition}[theorem]{Definition}

\begin{document}
\newcommand{\be}{\begin{equation}}
\newcommand{\ee}{\end{equation}}
\newcommand{\bs}{\begin{split}}
\newcommand{\es}{\end{split}}
\newcommand{\bt}{\begin{theorem}}
\newcommand{\et}{\end{theorem}}
\newcommand{\bd}{\begin{definition}}
\newcommand{\ed}{\end{definition}}
\newcommand{\bp}{\begin{proposition}}
\newcommand{\ep}{\end{proposition}}
\newcommand{\bl}{\begin{lemma}}
\newcommand{\el}{\end{lemma}}
\newcommand{\bc}{\begin{corollary}}
\newcommand{\ec}{\end{corollary}}
\newcommand{\bcon}{\begin{conjecture}}
\newcommand{\econ}{\end{conjecture}}
\newcommand{\la}{\label}
\newcommand{\A}{{\rm A}}
\newcommand{\rP}{{\rm P}}
\newcommand{\Z}{{\Bbb Z}}
\newcommand{\R}{{\Bbb R}}
\newcommand{\Q}{{\Bbb Q}}
\newcommand{\C}{{\Bbb C}}
\newcommand{\hra}{\hookrightarrow}
\newcommand{\lra}{\longrightarrow}
\newcommand{\lms}{\longmapsto}
\newcommand{\ra}{\rightarrow}
\newcommand{\sgn}{\rm sgn}
\newcommand{\Si}{\Sigma}
\begin{titlepage}
\title{Quantum   polylogarithms}   
\author{Alexander B. Goncharov}
\end{titlepage}

\setcounter{tocdepth}{3}

\date{}

\maketitle

\tableofcontents

 \begin{abstract} 
 
 Multiple polylogarithms  
are periods of variations of mixed Tate motives. Conjecturally, they deliver all such periods. 
 
 We introduce deformations of multiple polylogarithms depending on a parameter $\hbar\in \C$. 
 We call them {\it quantum  polylogarithms}.  Their asymptotic expansion as $\hbar \to 0$ recovers  multiple polylogarithms. 
 The quantum dilogarithm  was studied by Barnes in the XIX century.  
 Its exponent appears in  many areas of Mathematics and Physics. 
  
 Quantum polylogarithms 
 satisfy a holonomic systems of    
 modular difference equations with coefficients in  variations of mixed Hodge-Tate structures of motivic origin. 

 If $\hbar \in \Q$, the quantum polylogarithms  can be expressed via multiple polylogarithms.  
 
 However if $\hbar \not = \Q$, quantum polylogarithms  are not periods of variations of mixed motives, i.e. 
  they can not be  given by integrals of rational differential forms on algebraic varieties. 
 Instead, quantum polylogarithms are integrals of differential forms built  from both rational functions and 
 exponentials of rational functions.   We call them {\it rational exponential integrals}. 
 \vskip 1mm
 
 We suggest that quantum polylogarithms reflect a very general phenomenon: 
 \vskip 1mm
 {\it Periods of variations of mixed   motives should have  quantum deformations}.

 \end{abstract}

\section{Introduction}

\subsection{The dilogarithm and the quantum dilogarithm}
The dilogarithm function is defined by the following power series:
$$
{\rm Li}_2(z):= \sum_{k>0}\frac{z^k}{k^2}.
$$
It can be continued analytically to a function on a cover of ${\Bbb C}{\Bbb P}^1 - \{0, 1, \infty\}$ via the integral 
$$
{\rm Li}_2(z):= -\int_0^z\log(1-t)\frac{dt}{t}.
$$
 One of the key features of the dilogarithm is that it satisfies Abel's five term relation. \vskip 2mm

The dilogarithm power series  
admit a  $q-$deformation:  
$$
{\rm Li}_{1, 1}(z; q):= \sum_{k=1}^{\infty}\frac{z^k}{(q^{k} - q^{-k})  k}.
$$
Its exponent is the inverse of the Pochhammer symbol. Precisely, set 
$$
\Psi_q(z) := \frac{1}{(1+qz)(1+q^3z)(1+q^5z) \cdot  ~\ldots} .
$$
Then one has 
\be \la{PSI1I}
\log \Psi_{q}(z) = -{\rm Li}_{1, 1}(-z; q).
\ee
Setting $q=e^{\pi i \hbar}$, and letting $\hbar \to 0$, we get the asymptotic expansion 
\be \la{PhiPsi}
\log \Psi_{q}(z) \sim_{\hbar \to 0} ~ -\frac{{\rm Li}_2(-z)}{2\pi i \hbar}.
\ee
It satisfies the quantum pentagon relation,  discovered by Faddeev and Kashaev \cite{FK}. Namely, consider  variables $X, Y$ satisfying 
the relation $XY = q^2YX$. Then we have the identity  of $q-$commutative power series in $X, Y$:
$$
\Psi_q(X) \Psi_q(Y) = \Psi_q(Y) \Psi_q(X+Y)\Psi_q(X), 
$$
Its quasiclassical limit recovers  Abel's five term relation for the dilogarithm. 
\vskip 2mm

The power series  $\Psi_{q}(z) $  converge only  if   $|q|<1$. Remarkably, the  quotient 
\be \la {psi112}
{\Phi}_\hbar(\omega) = \frac{\Psi_q(e^\omega)}{\Psi_{q^\ast}(e^{\omega/\hbar})}, ~~~~q = e^{i\pi \hbar}, ~~q^\ast = e^{-i\pi /\hbar}, \ \ {\rm Im}\hbar >0
\ee
has excellent analytic properties. In particular, it is a  meromorphic function in $\omega$, depending on a complex parameter $\hbar$. 
To see this, recall the  integral  introduced and studied  by Barnes 
\cite{Ba}:
\be \la{in112}
\begin{split}
&{\cal F}^\hbar(\omega):= 
\int_{{\R+i0}}  
\frac{e^{-ip\omega}}{{\mathfrak{sh}}  (\pi p){\mathfrak{sh}}  (\pi \hbar p)}
\frac{d p}{p}, \ \ \ \ \ \  \ \ \ \ \ \ \mathfrak{sh}(p) =  e^p-e^{-p}.\\
\end{split}
\ee
The integration contour   $\R+i0$  is the limit of the contour $\R+i\varepsilon$, $\varepsilon >0$ 
 when $\varepsilon \to 0$.  The integral is well defined for any complex values of  $\hbar$,  convergent for ${\rm Im}(\omega) < \pi(1+ |{\rm Re}(\hbar)|)$, 
 and satisfies difference relations under the shift of  $\omega$ by $2\pi i$ and $2\pi i \hbar$
 \be
 \begin{split}
& {\cal F}^\hbar(\omega+2\pi i)  \ = \ {\cal F}^\hbar(\omega)  - \log (1+e^{i\pi \hbar}e^\omega) \\
& {\cal F}^\hbar(\omega+2\pi i\hbar) = {\cal F}^\hbar(\omega)  - \log (1+e^{i\pi/ \hbar}e^{\omega/\hbar})\\
 \end{split}
 \ee
 which allow to extend it   to a multivalued analytic  function in $\omega \in \C$. 
 Finally, one has   
\be
\Phi_\hbar(\omega)= {\rm exp}(-{\cal F}^\hbar(\omega)).
\ee
 We call the function ${\cal F}^\hbar(\omega)$ the {\it quantum dilogarithm},  although  this name is often used for  its exponent. The quantum dilogarithm  and its relatives appear in Statistical Physics \cite{Bax},  Liouville theory\cite{DO}, \cite{ZZ}, quantum groups \cite{F},  
  quantum higher Teichmuller theory \cite{K}, \cite{CF},  \cite{FKV}, \cite{FG1} -  \cite{FG2}, \cite{GS}, 
  quantization of cluster varieties \cite{FG3}, 
  and many other areas  \cite{V}. 
The {quantum dilogarithm} satisfies the   quantum pentagon relation, which plays an important role in its applications.  Namely, the function $\Phi_\hbar(\omega)$ is well defined and unitary on the real line. Consider the  unitary operator $K $ in  
$L_2(\R)$ given by the multiplication by $\Phi_\hbar(w)$, followed by the Fourier transform:
 $$
K: f(w) \lra \int_\R f(w) \Phi_\hbar(w) e^{-\frac{i\omega \xi}{2 \pi \hbar}} dt.
$$
Then $
K^5 = c \cdot {\rm Id}$, where  $c$ is a constant,  $|c|=1$. Its quasiclassical limit delivers the five term relation for the dilogarithm. 
The remarkable analytical properties of the function ${\cal F}^\hbar(\omega)$ together with the quantum deformation of the five term relation  convince  that it is the  natural quantum deformation of the dilogarithm. 
Let us now look at generalizations of the dilogarithm function. 

\subsection{Multiple polylogarithms} 
Recall the classical polylogarithm series:
$$
{\rm Li}_n(z):= \sum_{k>0}\frac{z^k}{k^n}.
$$
They make sense for any integer $n$. For $n\leq 0$ they are rational functions in $z$. 
For $n>0$ they are convergent for $|z|<1$, and 
admit an analytic continuation to a cover of  
$\C{\Bbb P}^1-\{0, 1, \infty\}$. 
\vskip 2mm
 
{\it Multiple polylogarithms}  are defined by the power series expansion \cite{G94}: 
\begin{equation} \label{5}
{\rm Li}_{n_{1},\ldots ,n_{m}}(z_{1},\ldots ,z_{m})   
:= \quad 
\sum_{0 < k_{1} < k_{2} < ... < k_{m} } \frac{z_{1}^{k_{1}}z_{2}^{k_{2}}
\ldots z_{m}^{k_{m}}}{k_{1}^{n_{1}}k_{2}^{n_{2}}\ldots k_{m}^{n_{m}}},
\end{equation}
which are convergent if $|z_{i}| < 1$. Here $m$ is the {\it depth}, and $|{\bf n}|:= n_1+\ldots +n_m$ is the {\it weight}. 

\vskip 3mm
Multiple polylogarithms admit an iterated integral presentation, 
which allows to continue them analytically. Namely,    given   meromorphic $1-$forms $\omega_i(t)$ on $\C$, and a path $\gamma:[0,1]\lra \C$  not intersecting their poles,  
recall the  iterated integrals on the line:
$$
\int_{\gamma}\omega_1(t) \circ \ldots \circ 
\omega_m(t) :=  
 \int_{0 \leq t_1 \leq \ldots \leq t_m\leq 1}\gamma^*\omega_1(t_1) \wedge  \ldots \wedge  
\gamma^*\omega_m(t_m). 
$$
To  present multiple polylogarithms by iterated integrals on the line,   set 
\begin{equation} \label{3??}
\begin{split}
&{\rm I}_{n_{1},\ldots ,n_{m}}(0; z_{1}, \ldots , z_{m}; z_{m+1}) : = 
 \int_{0}^{z_{m+1}} \underbrace {\frac{dt}{z_{1}-t} \circ \frac{dt}{t} \circ \ldots  
\frac{dt}{t}}_{ \mbox{$n_1$ differentials}} \circ 
    \ldots   \circ 
\underbrace {\frac{dt}{z_{m}-t} \circ \frac{dt}{t} \circ \ldots \circ \frac{dt}{t}}_
{\mbox {$n_m$ differentials}}.\\
\end{split}
\end{equation}
Here the iterated integral is over a path from $0$ to $z_{m+1}$. It is a multivalued analytic function on  the space of collections of distinct points $(0, z_1, ..., z_{m+1})$ in $\C$. 
Then by \cite[Theorem 2.2]{G01}:\footnote{
The {\it positive locus} ${\cal M}^+_{0, n+3}$ of the moduli space ${\cal M}_{0, n+3}$ the set of  ordered configurations of points  on ${\R}{\Bbb P}^{1}$  modulo the diagonal ${\rm PGL}_2(\R)$-action,  whose order is compatible with one of the   circle orientations.  
 We have 
 \be
(\infty,   z_1,    ..., z_m, 1, 0)\in {\cal M}^+_{0, n+3}\quad \mbox{if and only if} \quad z_1>  \ldots > z_m  >  1.
\ee
 Reversing the order  preserves the positive locus. Iterated integral (\ref{3??}) has a natural branch on the positive 
locus, 
  provided by the path 
$\gamma = [0, 1]$. 
The positivity just means  that arguments of  multiple polylogarithm series  (\ref{5*??})   are positive numbers smaller than $1$,   so the series are convergent. 
 }
\begin{equation} \label{5*??}
 {\rm I}_{n_{1},\ldots ,n_{m}}(0; z_{1}, \ldots , z_{m}; z_{m+1})~  = ~   {\rm Li}_{n_{1},\ldots ,n_{m}}
\Bigl(\frac{z_2}{z_1},\frac{z_3}{z_2},\ldots ,\frac{z_{m+1}}{z_m}\Bigr).  
\end{equation}

In Section \ref{SEC4.1} we establish  a new integral presentation for multiple polylogarithms:
\be \la{2.27.11.2n}
 \begin{split}
 & {\rm Li}_{{n_1}, \ldots , {n_m}}(e^{\omega_1-\omega_2}, e^{\omega_2-\omega_3}, \ldots , -e^{\omega_{m}}) =  \\
&i^{|{\bf n}|-m} \int_{(\R+i0)^m}
\frac{e^{-ip_1\omega_1}}{{\mathfrak{sh}}(\pi p_1)}\frac{dp_1}{p_1^{n_1}} \wedge \ldots \wedge 
\frac{e^{-ip_m\omega_m}}{{\mathfrak{sh}}(\pi p_m)}
\frac{dp_m}{(p_1+\ldots +p_m)^{n_m}}. \\
\end{split}
\ee
So we have three different presentations of  multiple polylogarithms: as power series (\ref{3??}), via iterated integrals (\ref{5*??}),  and 
 using integral presentation (\ref{2.27.11.2n}). We note  that the two  integrals (\ref{3??}) and (\ref{2.27.11.2n}) are related via power series (\ref{5}) rather than directly. 

  Iterated integral presentation (\ref{5*??}) shows that multiple polylogarithms are periods  of mixed Tate motives. Conjecturally, 
 they provide  all such periods \cite[Conjecture 17]{G94}. 
 \vskip 1mm
 
 We  introduce a deformation
 of    multiple polylogarithms, called  {\it quantum   polylogarithms}, depending on   a  
parameter $\hbar\in \C$. Using 
  integral presentation (\ref{2.27.11.2n}) we see  that its  asymptotic expansion  at $\hbar \to 0$  recovers multiple polylogarithms. 

Quantum polylogarithms   provide quantum deformation of all periods of mixed Tate motives.  

I suggest  that  periods of any variations of mixed motives admit quantum deformation.

\subsection{Quantum   polylogarithms.}   
Let us introduce first   the {\it kernel function}. 

\bd
Let $a, b$ be a pair of non-negative integers. The kernel function is defined by
 \be \la{EX1}
K^\hbar_{a,b}(p; \omega):= \frac{ e^{-ip\omega}}{{\mathfrak{sh}}^a(\pi p) {\mathfrak{sh}}^b(\pi \hbar p)}.
\ee
\ed

 The {\it depth $m$} quantum polylogarithms  are functions $ {\mathcal{F}}^\hbar_{{\bf a}, {\bf b}, {\bf n}}(\omega_1, \ldots , \omega_m)$ in  $m$ complex variables $\omega_i$ 
 which depend on  a triple of $m$-tuples of   integers
$$
{\bf a}= (a_1, ..., a_m), \quad {\bf b}= (b_1, ..., b_m), \quad {\bf n}= (n_1, ..., n_m), \qquad a_i,b_i \in \Z_{\geq 0}, ~ n_i \in \Z.
$$

\bd \la{DDQP} The depth $m$   quantum    polylogarithms  
 are the integrals
\be \la{15f}
\begin{split}
\mathcal{F}^\hbar_{{\bf a}, {\bf b}, {\bf n}}(\omega_1, ..., \omega_m):=  &
 i^{|{\bf n}|-m} \int_{(\R+i0)^m} 
\bigwedge_{k=1}^mK^\hbar_{a_k, b_k}(p_k; \omega_k)
\frac{dp_k}{(p_1+\ldots +p_k)^{n_k}}. \\\end{split}
\ee
We define the weight of $\pi $ and $\omega_j$ to be   1, and the {\it weight} of the quantum polylogarithm 
 to be
\be \la{n}
|{\bf n}|:= n_1 + ... + n_m.
\ee
\ed

 The integral converges if 
 $
|{\rm Im} ~\omega_i |< \pi(a_i+ b_i |{\rm Re} (\hbar)|)$. 
It  extends  to a multivalued analytic function in  $(\omega_1, ..., \omega_m)\in \C^m$ using the difference relations (\ref{2.20.11.2}) for quantum polylogarithms.
 
  As an analytic  function of $\hbar$, the 
 ${\cal F}_{ {\bf a}, {\bf b}, {\bf n}}^\hbar(\omega_1, ..., \omega_m)$ extends  
 to the  complex plane with the negative real axis $\hbar<0$ removed.  
 The integral converges for any $n_i\in \R$ if  one of the integers $a_i, b_i$ is positive.  
Here are two examples.

\begin{enumerate}
\item 
The depth one quantum polylogarithms are   the following 
integrals 
$$
{\mathcal{F}}_{a,b,n}^\hbar(\omega):= 
i^{n-1}\cdot \int_{{\R+i0}}  
\frac{e^{-ip\omega}}{{\mathfrak{s}}\mathfrak{h}^a (\pi p){\mathfrak{sh}}^b (\pi \hbar p)}
\cdot \frac{dp}{p^{n}}, \qquad a, b \geq 0.
$$

\item 

The depth two quantum   polylogarithms  are given by
\be
\begin{split}
&{{\cal F}}^\hbar_{{\bf a}, {\bf b}, {\bf n}}(\omega_1, \omega_2):=  
 \\ & i^{|{\bf n}|-2}\cdot  \int_{({\R+i0})^2}  
\frac{e^{-ip_1\omega_1}}{{\mathfrak{sh}}^{a_1}(\pi p_1){\mathfrak{sh}}^{b_1}(\pi \hbar p_1)}
\frac{e^{-ip_2\omega_2}}{{\mathfrak{sh}}^{a_2}(\pi p_2){\mathfrak{sh}}^{b_2}(\pi \hbar p_2)}\cdot 
\frac{dp_1}{p_1^{n_1}}\wedge \frac{dp_2}{(p_1+p_2)^{n_2}}. \\
\end{split}
\ee
\end{enumerate}

  \vskip 1mm

 \bt Quantum polylogarithms at the rational $\hbar \in \Q$  are periods of variations of mixed Tate motives of the same weight. \et
 
See the precise statement  in Theorem \ref{DRELaaa}. 
 In sharp contrast with this, quantum polylogarithms for irrational $\hbar\not \in \Q$ are {\bf not} periods of variations of mixed motives.
  \vskip 2mm
  
 By Theorem \ref{shuffle},  quantum   polylogarithms    
  satisfy   
shuffle product formulas, providing   an algebra structure.  
\vskip 1mm

Specialising $\omega_1=...=\omega_m=0$ and ${\bf a}={\bf b}=(1, ..., 1)$, we get 
an $\hbar-$deformation of the depth $m$ multiple $\zeta$-function. For example, for the depth $2$ we have
 \be
\begin{split}
&\zeta_\hbar(s_1, s_2):= 
 \int_{({\R+i0})^2}  
\frac{d p_1}{{\mathfrak{sh}} (\pi p_1){\mathfrak{sh}} (\pi \hbar p_1)p_1^{s-1}}  
\frac{d p_2}{{\mathfrak{sh}} (\pi p_2){\mathfrak{sh}} (\pi \hbar p_2) (p_1+p_2)^{s_2-1}} . \\
\end{split}
\ee 
The analytic continuation of these functions is obtained the same way as in \cite[Theorem 2.25]{G01}. 
When $s_i=n_i$ are positive integers, we get an $\hbar-$deformation of Euler's multiple $\zeta$-values. \vskip 1mm

  Interesting $q-$deformations of the multiple $\zeta-$values were considered by Okounkov \cite{O}.

\subsection{$q-$deformations  of multiple polylogarithms} \la{CQMP}

{\it Convention}. We denote by latin letters $z_i$  coordinates of the points on the projective line, and by the greek letters $\omega_i$ the relevant logarithmic coordinates: $z_i=e^{\omega_i}$.
\vskip 1mm

 Multiple polylogarithm power series (\ref{5}) have a $q-$deformation:

\bd \la{MQP} Let  ${\bf a}, {\bf n}\in \Z_{\geq 0}^m$. Multiple $q-$polylogarithms are   power series in $z_1, ..., z_m$:
\be \la{5q}
\begin{split}
&{\rm Li}_{{\bf a}, {\bf n}}(z_1, ..., z_m; q):= \\
&\sum_{ k_1,  ...,k_m>0}^{\infty}
\frac{z_1^{k_1}z_2^{k_2}\ldots z_m^{k_m}}{[k_1]_q^{a_1} [k_2]_q^{a_2} \ldots 
[k_m]_q^{a_m}\cdot k_1^{n_1} (k_1+k_2)^{n_2} \ldots (k_1+...+k_m)^{n_m}}, \quad [n]_q:= q^n-q^{-n}.\\
\end{split}
\ee
\ed

The weight of the multiple $q-$polylogarithm series (\ref{5q})  is defined to be ${\bf n}$, just as in (\ref{n}). 

{\it Example}.  The weights of  the $q-$dilogarithm 
${\rm Li}_{1,1}(z; q)$ and the  quantum dilogarithm ${\cal F}^\hbar(\omega) $ are $1$.  Since the weight of $2\pi \hbar$ is   $1$, the weight of the function  in (\ref{PhiPsi}) is 
$1$, consistently with (\ref{PSI1I}).  
\vskip 1mm

Note that unlike in series (\ref{5}), the summation in (\ref{5q}) is over the octant $k_1, ..., k_m>0$. 
Therefore there are   no shuffle product formulas for the series ${\rm L}_{{\bf a}, {\bf n}}(z_1, ..., z_m; q)$.  

Power series (\ref{5q})    are not defined when $q$ is a root of unity. 
   In Section \ref{sec3.2} we complement them by {\it companion series}. Their appropriate  sums  coincide with the  quantum polylogarithm  integrals, generalising the logarithm of relation (\ref{psi112}), 
 thus converging for any $\hbar$.  
 
\subsection{Connections between quantum and  multiple polylogarithms}

The  weight is compatible with  major operations with quantum polylogarithms:

\begin{itemize}

\item The partial derivatives   decrease the weight by $1$, see (\ref{2.20.11.1}). 

\item The difference relations preserve the weight, see (\ref{2.20.11.2}).

\end{itemize} 
 Passing from quantum polylogarithms to their exponents destroys the  weights. In particular, this is why we call the function ${\cal F}^\hbar(\omega)$, rather than  its exponent $\Phi_\hbar(\omega)$, the
   quantum dilogarithm. \\
   
   Quantum polylogarithms are related to  multiple polylogarithms  in several  ways:

\begin{enumerate}  
 \item {\it Via the asymptotic expansion as $\hbar \to 0$}. The weight $|{\bf n}|$ quantum polylogarithms  
have an  asymptotic expansion as $\hbar \to 0$ of the following shape, see Theorem \ref{asexp}: 
$$
  \sum_k (2 \pi i \hbar)^{-k} \times ~~\mbox{sums of multiple polylogarithms of the weight $|{\bf n}|+k$}. 
  $$
  All terms of the asymptotic expansion have the same weight $|{\bf n}|$. 
  
 \item {\it Via the $\hbar=1$ specialization}.  We prove in  Theorem \ref{h=1}:
    \bt \la{h=1a} 
The 
${{\cal F}}^1_{{\bf a}, {\bf b}, {\bf n}}(\omega_1, ..., \omega_m)$ is the product of  the function ${\bf P}_{{\bf a}, {\bf b}}\Bigl(\frac{\omega_1}{2\pi}. ..., \frac{\omega_m}{2\pi}\Bigr)$, where 
${\bf P}_{{\bf a}, {\bf b}} $ is a polynomial with coefficients in $\Q$, and 
the multiple polylogarithm ${\rm Li}_{n_1, ..., n_m}$:
\be \la{h=1e}
{{\cal F}}^1_{{\bf a}, {\bf b}, {\bf n}}(\omega_1, ..., \omega_m) = {\bf P}_{{\bf a}, {\bf b}}\Bigl(\frac{\omega_1}{2\pi}. ..., \frac{\omega_m}{2\pi}\Bigr) 
\cdot {\rm Li}_{n_1, ..., n_m}\Bigl((-1)^{a_1+b_1}e^{\omega_1}, ..., (-1)^{a_m+b_m}e^{\omega_m}\Bigr).\ee
\et
Since the weight of $\frac{\omega}{2\pi}$ is zero, both parts of the equality have the same weight $|{\bf n}|$. 

The right hand side of (\ref{h=1e}) is a period of a variation of mixed Tate motives on $(\C^\times)^m$:
$$
 (\ref{h=1e})=  {\bf P}_{{\bf a}, {\bf b}}\Bigl(\frac{\log z_1}{2\pi}. ..., \frac{\log z_m}{2\pi}\Bigr) 
\cdot {\rm Li}_{n_1, ..., n_m}\Bigl((-1)^{a_1+b_1}z_1, ..., (-1)^{a_m+b_m}z_m\Bigr).
$$

   \item {\it Via   distribution relations}. Given a pair of coprime integers $r, s$, Theorem \ref{DREL}  relates quantum polylogarithms at  $\frac{r}{s}\hbar$ to a sum of similar quantum polylogarithms at $r \hbar$. Precisely:
   
    \bt \la{DRELaa} One has 
 distribution relations:
\be \la{DRELa}
\begin{split}
&r^{|{\bf n}|-m}{\cal F}_{{\bf a}, {\bf b}, {\bf n}}^{\frac{r}{s}\hbar}(r\omega_1, ..., r\omega_m) = \\
&\sum\limits_{\alpha_j=\frac{1-r}{2}}^{\frac{r-1}{2}} 
\sum\limits_{\beta_j=\frac{1-s}{2}}^{\frac{s-1}{2}} {\cal F}_{{\bf a}, {\bf b}, {\bf n}}^\hbar(\ldots, \omega_k+\frac{2\pi i}{r}\sum_{i=1}^{a_k}\alpha_j
+ \frac{2\pi i\hbar }{s}\sum_{j=1}^{b_k}\beta_j , \ldots).\\
\end{split}
\ee 
Here the sum is   over half-integers $\alpha_j, \beta_j$ if the summation limits are 
half-integers.  
 \et

Combining Theorems \ref{h=1a} and  \ref{DRELaa}, we express quantum polylogarithms at the rational   $\hbar\in \Q$ via sums of the multiple polylogarithms of the same weight:
 
    \bt \la{DRELaaa} Quantum polylogarithm functions at $\hbar \in \Q$ are  periods of variations of mixed Tate motives, provided   by  multiple polylogarithms. Precisely, set 
$$
\omega_k'= \omega_k +\frac{2\pi i}{r}\sum_{i=1}^{a_k}\alpha_j
+ \frac{2\pi i\hbar }{s}\sum_{j=1}^{b_k}\beta_j.$$
Then, using the notation of Theorem \ref{h=1a}, we have 
 \be \la{DRELa}
\begin{split}
&r^{|{\bf n}|-m}{\cal F}_{{\bf a}, {\bf b}, {\bf n}}^{\frac{r}{s}}(r\omega_1, ..., r\omega_m)  \\
& = \sum\limits_{\alpha_j=\frac{1-r}{2}}^{\frac{r-1}{2}} 
\sum\limits_{\beta_j=\frac{1-s}{2}}^{\frac{s-1}{2}} {\bf P}_{{\bf a}, {\bf b}}\Bigl(\frac{\omega'_1}{2\pi}. ..., \frac{\omega'_m}{2\pi}\Bigr) 
 {\rm Li}_{n_1, ..., n_m}(e^{\omega_1'}, ..., e^{\omega_m'}).\\
\end{split}
\ee 

 \et 

\item {\it Via    companion $q-$polylogarithm series}. Recall that we have 
\be \la{po}
\begin{split}
{\cal F}^\hbar(\omega) ~\stackrel{(\ref{psi112})}{=}  ~&\log \Psi_q(e^\omega) - \log \Psi_{q^\ast}(e^{\omega/\hbar})\\
 ~\stackrel{(\ref{PSI1I})}{=} ~&-{\rm Li}_{1,1}(-e^\omega; q) +  {\rm Li}_{1,1}(-e^{\omega/\hbar}; q^*). \\
\end{split}
\ee 
The depth $m$ quantum polylogarithms are sums of $2^m$ {\it companion $q-$polylogarithm series}  of the same weight, 
 generalizing   formula (\ref{po}).

\end{enumerate}

So the weight of quantum polylogarithms is compatible, in several different ways,   with the weight of multiple polylogarithms.
This is quite remarkable since the weight of multiple polylogarithms have  a deep algebraic geometric origin, while quantum polylogarithms  live outside of the traditional  Algebraic Geometry.

\vskip 2mm
{\bf Acknowledgement}. This work was supported by the NSF grant DMS-2153059. I am grateful to IHES for hospitality and support in various stages of the preparation of the paper.

\section{A new integral presentation  for multiple polylogarithms} \la{SEC4.1}


\bt \la{3.7.11.1}
  Assume that $|{\rm Im}~w_i| < \pi$ and 
${\rm Re}~w_i <0.$ Then one has
\be \la{2.27.11.2a}
\begin{split}
& {\rm Li}_{{n_1}, \ldots , {n_m}}(e^{w_1},    \ldots,  e^{w_{m-1}}, -e^{w_m})=\\
&i^{|{\bf n}|-m}  \int_{(\R+i0)^m} 
\frac{e^{-ip_1w_1}}{{\mathfrak{sh}}(\pi p_1)}\frac{dp_1}{p_1^{n_1}} \wedge \ldots \wedge 
\frac{e^{-i(p_1+...+p_m)w_m}}{{\mathfrak{sh}}(\pi p_m)}
\frac{dp_m}{(p_1+\ldots +p_m)^{n_m}}. \\
\end{split}
\ee
\et

\begin{proof}  

Consider 
 the integral over   
$(\Omega_N+i\varepsilon)^m$, where $N>0$ is an integer,   
$\Omega_N$ is the square in the upper half plane with the base 
$[-N, N]$, and $(\Omega_N+i\varepsilon)$   
its shift  by a small $\varepsilon>0$. 

We claim that this integral is convergent, and the integrals over any but 
the bottom side  decay exponentially as $N \to \infty$.  
 Indeed,  we have
\be
\begin{split}
&|e^{-ipw}| = e^{{\rm Im}(p){\rm Re}(w)+{\rm Re}(p){\rm Im}(w)}, \\
& |e^{-ipw}/{\mathfrak{sh}}(\pi p)| \sim_{p\to \pm \infty}  e^{{\rm Im}(p){\rm Re}(w)+{\rm Re}(p)({\rm Im}(w)\mp \pi)}.\\
\end{split}
\ee 
So  $e^{-ipw}/{\mathfrak{sh}}(\pi p)$ decays exponentially on the left and right sides as $N \to \infty$ 
since ${\rm Re}~w<0$. 
The integral over the top side of $(\Omega_N+i\varepsilon)^m$ decays 
exponentially since   $|{\rm Im}~w| < \pi$, and therefore 
$$
{\rm Re}(p)({\rm Im}(w)\mp \pi)\lra -\infty~~ \mbox{if} ~~ p\lra \pm \infty.
$$
 Indeed, we have either ${\rm Re}(p) \to \infty$ $\&$ ${\rm Im}(w)- \pi<0$, or  ${\rm Re}(p) \to -\infty$ $\&$ ${\rm Im}(w)+ \pi>0$. 
Finally, on our contour 
  $|p| > \varepsilon$, so the integral over the bottom side converges.   

Therefore we can calculate the integral  using the residue theorem. 
The residues are at the points 
$(p_1, ..., p_m)=(ik_1, ..., ik_m)$, where $k_1, ..., k_m>0$ are  integers.\footnote{See also the proof of Theorem \ref{h=1} where we explain  how the calculation of residues  in the more general set up at $p_1=ik_1, ..., p_m=ik_m$ where $k_1, ..., k_m>0$  
reduces to the calculation of the residues at $p_1=...=p_m=0$.} For example, in the depth $2$ case the contribution of the residue at the point $(p_1, p_2)= (ik_1, ik_2)$  is equal 
to\footnote{The factor $i^{|{\bf n}|}$ in (\ref{2.27.11.2a}) cancells with the factor $i^{|{\bf n}|}$ from the denominator. This is how  the factor $i^{|{\bf n}|}$ in  (\ref{2.27.11.2a}), as well as  in Definition \ref{DDQP} of quantum polylogarithms, helps. Next, the factor $(2\pi i)^m$ from the Cauchy theorem cancell the  factor $i^{-m}$ in  (\ref{2.27.11.2a}) and Definition \ref{DDQP}.}
\be \la{rescal}
(2\pi i)^2\cdot \frac{i^{-2}}{(2\pi)^2}\frac{(-1)^{k_1+k_2}e^{k_1w_1}e^{w_2(k_1+k_2)}}{k_1^{n_1}(k_1+k_2)^{n_2}} =   \frac{ e^{k_1w_1}(-e^{w_2})^{k_1+k_2}}{k_1^{n_1}(k_1+k_2)^{n_2}}.
\ee
The  $(-1)^{k_1+k_2}$ amounts to the fact that $\mathfrak{sh}(\pi p_j)$ is multiplied by $(-1)^k$ after the shift by $i\pi k$. Then the sum   $\sum_{k_1>0, k_2>0}$  delivers  
${\rm Li}_{{n_1}, {n_2}}(e^{w_1}, -e^{w_2})$. 
The series are convergent since ${\rm Re}(w_i)<0$. 

The argument in the depth $m$ case is the similar: the residue at $(ik_1, ..., ik_m)$ is 
\be \la{rescala}
\begin{split}
 (2\pi i)^m\cdot \frac{i^{-m}}{(2\pi)^m}&\frac{ (-1)^{k_1+...+k_m}e^{k_1w_1}\ldots e^{(k_1+... + k_m)w_m}}{k_1^{n_1}\ldots (k_1+...+k_m)^{n_m}} \\
 =&\frac{ e^{k_1w_1} e^{(k_1+k_{2})w_{2}} \ldots(-e^{w_m})^{k_1+... + k_m}}{k_1^{n_1}(k_1+k_{2})^{n_{2}}\ldots (k_1+...+k_m)^{n_m}}.\\
\end{split}
\ee
Theorem is proved.
\end{proof}

Let $\gamma = [0, 1]$. Denote by ${\rm I}_{{n_1}, \ldots , {n_m}}^{(\gamma)}(z_1, \ldots, z_m)$ the iterated integral defined using this path. 
\bc  \la{3.2} Assume that $|{\rm Im}~\omega_i| < \pi$ and 
\be \la{2.25.11.1}
 {\rm Re}~\omega_1 < {\rm Re}~\omega_2< \ldots <{\rm Re}~\omega_m<0.
\ee 
Then one has
\be \la{2.27.11.2}
 \begin{split}
 &{\rm I}^{(\gamma)}_{{n_1}, \ldots , {n_m}}(0; e^{-\omega_1}, \ldots , e^{-\omega_{m}};  -1)\stackrel{(\ref{5*??})}{=}\\
 & {\rm Li}_{{n_1}, \ldots , {n_m}}(e^{\omega_1-\omega_2}, e^{\omega_2-\omega_3}, \ldots , -e^{\omega_{m}}) =  \\
&i^{|{\bf n}|-m} \int_{(\R+i0)^m}
\frac{e^{-ip_1\omega_1}}{{\mathfrak{sh}}(\pi p_1)}\frac{dp_1}{p_1^{n_1}} \wedge \ldots \wedge 
\frac{e^{-ip_m\omega_m}}{{\mathfrak{sh}}(\pi p_m)}
\frac{dp_m}{(p_1+\ldots +p_m)^{n_m}}. \\
\end{split}
\ee
\ec 

\begin{proof} 
The first equality is the basic equality (\ref{5*??}). The second gives  the new integral presentation (\ref{2.27.11.2n}) for multiple polylogarithms, and  follows immediately from     (\ref{2.27.11.2a}). Indeed, 
the integrands in (\ref{2.27.11.2a}) and (\ref{2.27.11.2}) differ only by the exponentials, which match thanks to the identity
\be \la{popi1}
\begin{split}
& e^{-ip_1\omega_1} e^{-ip_2\omega_2} \ldots  e^{-ip_m\omega_m} = e^{-ip_1(\omega_1-\omega_2)}e^{-i(p_1+p_2)(\omega_2-\omega_3)} \ldots e^{-i(p_1+...+p_m)\omega_m}.\\
\end{split}
\ee
Note that the latter is equivalent to the identity 
\be \la{popi2a}
\begin{split}
&p_1 \omega_1 + p_2\omega_2 + \ldots + p_m\omega_m= p_1(\omega_1-\omega_2) + (p_1+p_2) (\omega_2-\omega_3) + ... + (p_1+... +p_{m})\omega_m. \\
\end{split}
\ee
Note that the condition on the  $\omega_i$ in Corollary \ref{3.2} is  equivalent to  the one  in Theorem \ref{3.7.11.1}. \end{proof}

So we get new  
  integral presentations for both multiple polylogarithms and iterated integrals. For example, for  the depth $m=2$ they look as follows:
  \be
\begin{split}
&{\rm Li}_{{n_1}, {n_2}}(e^{w_1}, -e^{w_2})= i^{|{\bf n}|-2}\int_{({\R+i 0})^2}  
\frac{e^{-ip_1w_1}}{{\mathfrak{sh}}(\pi p_1)}
\frac{e^{-i(p_1+p_2)w_2}}{{\mathfrak{sh}} (\pi p_2)} 
\frac{dp_1}{p_1^{n_1}}\wedge \frac{dp_2}{(p_1+p_2)^{n_2}}. \\
&{\rm I}_{{n_1}, {n_2}}(0; e^{-\omega_1}, e^{-\omega_2}; -1)= i^{|{\bf n}|-2}\int_{(\R+i0)^2}  
\frac{e^{-ip_1\omega_1}}{{\mathfrak{sh}}(\pi p_1)}
\frac{e^{-ip_2\omega_2}}{{\mathfrak{sh}}(\pi p_2)}
\frac{dp_1}{p_1^{n_1}}\wedge \frac{dp_2}{(p_1+p_2)^{n_2}}. \\
\end{split}
\ee

  \paragraph{\it Remark.}  
  
Condition (\ref{2.25.11.1}) for real   $\omega_i$ just means that  $e^{-\omega_1} >   \ldots > e^{-\omega_m}>1$. 
So   the arguments of the iterated integral in (\ref{2.27.11.2}) form a positive configuration of $m+3$ points. 

 \section{Properties of quantum   polylogarithms}

\subsubsection{Difference relations} Recall  the kernel function 
 \be \la{EX1}
K^\hbar_{a,b}(p; \omega):= \frac{ e^{-ip\omega}}{{\mathfrak{sh}}^a(\pi p) {\mathfrak{sh}}^b(\pi \hbar p)}.
\ee
Let $\Delta^{(\omega)}_{a}$ be 
  difference operators in the variable $\omega$:
\be \la{difO}
\begin{split}
&\Delta^{(\omega)}_{a}f(\omega):= f(\omega+a) - f(\omega-a),\\
 \end{split}
\ee
The kernel function satisfies two difference equations in $\omega$:
\be \la{EX2}
\begin{split}
&\Delta^{(\omega)}_{i\pi}K^\hbar_{a,b}(p; \omega) ~=~
K^\hbar_{a-1,b}(p; \omega), \\
&\Delta^{(\omega)}_{i\pi \hbar}K^\hbar_{a,b}(p; \omega) ~=~
K^\hbar_{a,b-1}(p; \omega).  \\
\end{split}
\ee
Indeed,  one has
\be \la{dek}
\begin{split}
&\Delta^{(\omega)}_{i\pi \hbar}e^{-ip\omega}  ~=~ \mathfrak{sh}(\pi \hbar p)\cdot e^{-ip\omega} , \\
&\Delta^{(\omega)}_{i\pi }e^{-ip\omega}    ~= ~\mathfrak{sh}(\pi p)\cdot e^{-ip\omega} . \\
\end{split}
\ee 
Set   ${\bf 1}_k:= (0, ..., 0, 1, 0, ..., 0)$, where $1$ is on the $k$-th place. Then   difference relations (\ref{EX2}) for the kernel function  imply 
difference relations for quantum polylogarithms:\footnote{We  assume that ${\bf a} - {\bf 1}_k$ and ${\bf b} - {\bf 1}_k$ are still non-negative integers, to avoid convergence issues. }
\be \la{2.20.11.2}
\begin{split}
&\Delta^{(\omega_k)}_{i\pi}{{\cal F}}^\hbar_{{\bf a}, {\bf b}, {\bf n}}(\omega_1, ..., \omega_m) = 
{{\cal F}}^\hbar_{{\bf a} - {\bf 1}_k, {\bf b}, {\bf n}}(\omega_1, ..., \omega_m), \\
&\Delta^{(\omega_k)}_{i\pi\hbar}{{\cal F}}^\hbar_{{\bf a}, {\bf b}, {\bf n}}(\omega_1, ..., \omega_m) = 
{{\cal F}}^\hbar_{{\bf a} , {\bf b}- {\bf 1}_k, {\bf n}}(\omega_1, ..., \omega_m). \\
\end{split}
\ee

\subsubsection{The asymptotic expansion when $\hbar \to 0$} 

\bt \la{asexp} When $\hbar \to 0$,  the function 
 $
 {{\cal F}}^\hbar_{{\bf a}, {\bf b}, {\bf n}}(\omega_1, ..., \omega_m)
 $
has an asymptotic Laurent series expansion  in  $2\pi \hbar$,   
whose coefficients are sums of quantum  polylogarithms:
\be \la{36}
{{\cal F}}^\hbar_{{\bf a}, {\bf b}, {\bf n}}(\omega_1, ..., \omega_m) \sim_{\hbar \to 0} (2\pi \hbar)^{-|{\bf b}|}{{\cal F}}^\hbar_{{\bf a}, {\bf 0}, {\bf b}+ {\bf n}}(\omega_1, ..., \omega_m)+ \ldots
\ee
All terms of the asymptotic expansion have the same weight. 
\et

\begin{proof} The leading term of the $\hbar \to 0$ asymptotic expansion   of the kernel function (\ref{EX1})  is    
 \be \la{EX1a}
  \frac{ e^{-ip\omega}}{\mathfrak{sh}^a(\pi p){\mathfrak{sh}}^b(\pi \hbar p)} ~\sim_{\hbar \to 0}~ \frac{ 1}{(2 \pi \hbar)^b } \cdot  \frac{e^{-ip\omega}}{{\mathfrak{sh}}^a(\pi p)}\frac{1}{p^b} .
\ee 
This implies the claim about the leading term of the asymptotic expansion. For example, 
in  the depth two case 
the leading term of the expansion is   
$$
\frac{  i^{|{\bf n}|-2}}{(2\pi \hbar)^{|{\bf b}| }} \int_{(\R+i0)^2}   
\frac{e^{-ip_1\omega_1}}{{\mathfrak{sh}}^{a_1}(\pi p_1)}
\frac{e^{-ip_2\omega_2}}{{\mathfrak{sh}}^{a_2}(\pi p_2)}
\frac{dp_1}{p_1^{b_1+n_1}}\frac{dp_2}{p_2^{b_2}(p_1+p_2)^{n_2}}. 
$$
To get all terms of the asymptotic expansion, we write the Laurent series expanding $\mathfrak{sh}(\pi \hbar p_k)$ in $\pi \hbar p_k$. Since the weight of  $\pi \hbar p_k$ is zero, the weight conservation  is clear. 

To get the rest of the terms of the asymptotic expansion, we use   identity 
\be \la{pp+}
\frac{1}{p_2(p_1+p_2)}  = \frac{1}{p_1}\Bigl(\frac{1}{p_2} - \frac{1}{p_1+p_2}\Bigr).  
\ee
and proceed by the induction on $b_2+n_2$, till either $b_2=0$ or $n_2=0$. If $b_2=0$, we get the double polylogarithm. If $n_2=0$, we get a product of two depth one polylogarithms. 

The  case $m>2$ is similar. \end{proof}

\subsubsection{Distribution relations}

\bt \la{DREL} One has 
 distribution relations:
\be
\begin{split}
&r^{|{\bf n}|-m}{\cal F}_{{\bf a}, {\bf b}, {\bf n}}^{\frac{r}{s}\hbar}(r\omega_1, ..., r\omega_m) = \\
&\sum\limits_{\alpha_j=\frac{1-r}{2}}^{\frac{r-1}{2}} 
\sum\limits_{\beta_j=\frac{1-s}{2}}^{\frac{s-1}{2}} {\cal F}_{{\bf a}, {\bf b}, {\bf n}}^\hbar(\ldots, \omega_k+\frac{2\pi i}{r}\sum_{i=1}^{a_k}\alpha_j
+ \frac{2\pi i\hbar }{s}\sum_{j=1}^{b_k}\beta_j , \ldots).\\
\end{split}
\ee 
Equivalently, the  functions  ${\cal F}_{{\bf a}, {\bf b}, {\bf n}}^\hbar(\zeta_1, ..., \zeta_m)$ in the bottom line have the arguments 
$$
\zeta_k:= \omega_k+\frac{2\pi i}{r}\sum_{i=1}^{a_k}\alpha_j
+ \frac{2\pi i\hbar }{s}\sum_{j=1}^{b_k}\beta_j, \ \ \ \ \ \ k=1, ...,m.
$$
Here the sum is   over half-integers $\alpha_j, \beta_j$ if the summation limits are 
half-integers.  
 \et

 \begin{proof} Write the identity
$
  {\mathfrak{sh}}(rx) = {\mathfrak{sh}}( x) (e^{({r-1} )x} + e^{{(r-3)}x} + \ldots +e^{(1-r) x} )
 $ 
  as 
\be \la{idp}
\frac{1}{{\mathfrak{sh}}( x)}=  \frac{e^{({r-1} )x} + e^{{(r-3)}x} + \ldots +e^{(1-r) x}}{  {\mathfrak{sh}}(rx) }.
\ee
  Set $q=pr$ in the kernel function:
$$
K_{a,b}^{r/s}(p; r\omega):= \frac{ e^{-ipr\omega}}{{\mathfrak{sh}}^a(\pi p) {\mathfrak{sh}}^b(\pi \hbar p)} =  \frac{ e^{-iq\omega}}{{\mathfrak{sh}}^a(\pi q/r) {\mathfrak{sh}}^b( \pi q/s)}.
$$
Then using (\ref{idp}) 
we write this as 
\be
\begin{split}
&\frac{ e^{-iq\omega}\Bigl(e^{\frac{r-1}{r}\pi q} + e^{\frac{r-3}{r}\pi q} + \ldots +e^{\frac{1-r}{r} \pi q}\Bigr)^a
\Bigl(e^{\frac{s-1}{s}\pi q} + e^{\frac{s-3}{s}\pi q} + \ldots +e^{\frac{1-s}{s} \pi q}\Bigr)^b}{{\mathfrak{sh}}^a(\pi q) {\mathfrak{sh}}^b( \pi q)}\\
\end{split}
\ee The claim follows immediately from this by expanding the  products. 
\end{proof}

\subsubsection{The value at $\hbar=1$}  \la{secdistr}
Recall the basic quantum polylogarithm functions ${{\cal F}}^\hbar_{n_1, ..., n_m}(\omega_1, ..., \omega_m)$, see (\ref{153}).

  \bt \la{h=1} 
The function 
${{\cal F}}^1_{{\bf a}, {\bf b}, {\bf n}}(\omega_1, ..., \omega_m)$ is the product of a polynomial ${\bf P}_{{\bf a}, {\bf b}}\Bigl(\frac{\omega_1}{2\pi}. ..., \frac{\omega_m}{2\pi}\Bigr)$ with rational coefficients, and a
multiple polylogarithm of the same  depth $\&$ weight:
\be
{{\cal F}}^1_{{\bf a}, {\bf b}, {\bf n}}(\omega_1, ..., \omega_m) = {\bf P}_{{\bf a}, {\bf b}}\Bigl(\frac{\omega_1}{2\pi}. ..., \frac{\omega_m}{2\pi}\Bigr) 
\cdot {\rm Li}_{n_1, ..., n_m}\Bigl((-1)^{a_1+b_1}e^{\omega_1}, ..., (-1)^{a_m+b_m}e^{\omega_m}\Bigr)\ee
\et

\begin{proof}One has 
\be \la{098}
  \begin{split}
&{{\cal F}}^1_{{\bf a}, {\bf b}, {\bf n}}(\omega_1, ..., \omega_m) =  
i^{|{\bf n}|-m}  \int_{(\R+i0)^m} 
\prod_{j=1}^m\frac{e^{-ip_j\omega_j}}{{\mathfrak{sh}}^{a_j+b_j}(\pi p_j)}
\frac{dp_j}{(p_1+\ldots +p_j)^{n_j}}. \\
\end{split}
  \ee
Let us assume that ${\rm Re}(\omega_j)<0$. Then if ${\rm Im}(p_j) \to +\infty$, the exponential $e^{-ip_j\omega_j}$ decays fast. Next, when $|{\rm Re}(\omega_j)|< \pi(a_j+b_j)$, the integrand decays exponentially at $|p|\to \infty$. So assuming 
$
-\pi(a_j+b_j) <{\rm Re}(\omega_j)<0
$ 
we 
evaluate the integral as $(2\pi i)^m\times$  the sum over $k_1, ..., k_m>0$ of the residues at  
  $p_1=ik_1, ..., p_m=ik_m$. 
 \bl  Such a residue at  
  $p_1=ik_1, ..., p_m=ik_m$ is equal to
\be \la{765}
\begin{split}
&i^{-m}(-1)^{(a_1+b_1)k_1+ ... + (a_m+b_m)k_m}~{\rm Res}_{p_1=...=p_m=0} \Bigl(\prod_{j=1}^m\frac{e^{-ip_j\omega_j}dp}{{\mathfrak{sh}}^{a_j+b_j}(\pi p_j)}\Bigr)\cdot\frac{ e^{k_1\omega_1}\ldots e^{k_m\omega_m}}{k_1^{n_1} \ldots k_m^{n_m}}.\\
\end{split}
\ee
\el

\begin{proof}
Calculating the residue at $p_1=ik_1, ..., p_m=ik_m$ we use the following:

\begin{enumerate}

\item  The  shift by  $p_j \to p_j+i\pi$ for $j=1, ..., m$ results in the multiplication of the  denominator  $\prod_{j=1}^m{\mathfrak{sh}}^{a_j+b_j}(\pi p_j)$   by $(-1)^{|{\bf a}|+|{\bf b}|}$.
 
\item We expand the exponential at $p_j= ik_j+p'_j$ as 
$$
\prod_{j=1}^me^{-ip_j\omega_j} = \prod_{j=1}^me^{-ip'_j\omega_j} \cdot e^{k_1\omega_1}\ldots e^{k_m\omega_m}.
$$

\item  As $p_j \to ik_j$, we have
$$
 p_1^{-n_1} \cdot \ldots \cdot (p_1+...+p_m)^{-n_m} \lra i^{-|{\bf n}|}k_1^{-n_1} \cdot \ldots \cdot (k_1+...+k_m)^{-n_m}.
 $$ 
 \end{enumerate}
In particular we observe that the factor $i^{|{\bf n}|}$ in (\ref{098}) is cancelled with the one $i^{-|{\bf n}|}$ from (3).
 \end{proof}
 
 The left factor in (\ref{765}) does not depend on $k_j$. So  taking the sum over all $k_1, \dots, k_m$ we get 
\be
\begin{split}
& (2\pi)^m{\rm Res}_{p_1=...=p_m=0} \Bigl(\prod_{j=1}^m\frac{e^{-ip_j\omega_j}}{{\mathfrak{sh}}^{a_j+b_j}(\pi p_j)}\Bigr)\cdot 
{\rm Li}_{n_1, \ldots, n_m}
\Bigl((-1)^{a_1+b_1}e^{\omega_1}, \ldots,  (-1)^{a_k+b_k}e^{\omega_m}\Bigr).     
\\
\end{split}
\ee 
The residue on the left factorises into the product of the one variable residues:
$$
2\pi ~{\rm Res}_{p_j=0} \Bigl(\frac{e^{-ip_j\omega_j}dp_j}{{\mathfrak{sh}}^{a_j+b_j}(\pi p_j)} \Bigr)= 2\pi ~{\rm Res}_{p_j=0} \Bigl(\frac{1 + (-ip_j\omega_j) +   (-ip_j\omega_j)^{2} /2! + ...  }{ (\pi p_j)^{a_j+b_j}} (1 + ...)\Bigr).
$$
The latter  is 
equal to a sum of rational constants times 
$$
2\pi ~{\rm Res}_{p_j=0}\Bigl(\frac{ (-ip_j\omega_j)^{m}(\pi p_j)^ndp_k}{(\pi p_j)^{a_j+b_j}}\Bigr), \ \ \ \ \ \ m+n= a_j+b_j -1.
$$
Its weight is equal to zero since the weights of $\pi$ and $\omega_j$   are equal to $1$. 

Note also that $z_j=e^{\omega_j}$ has zero weight: indeed, we have $w_j = \log(z_j)$, so in the final answer we have polynomials in $z_j, \log (z_j)$.
 
\end{proof}

\subsubsection{The ${\cal I}-$variant of quantum   polylogarithms.}  \la{secty}
 
  We will need  the ${\cal I}-$variant of quantum   polylogarithms:
\be \la{SECTSS}
\begin{split}
&{{\cal I}}^\hbar_{{\bf a}, {\bf b}, {\bf n}}(w_1, ..., w_m):=
     i^{|{\bf n}|-m}  \int_{({\R+i0})^m} 
\prod_{k=1}^m\frac{e^{-i(p_1+... + p_k)w_k}}{{\mathfrak{sh}}^{a_k} (\pi p_2) \cdot {\mathfrak{sh}}^{b_k}( \pi \hbar p_2)}
\frac{dp_k}{(p_1+\ldots +p_k)^{n_k}}. \\
\end{split}
\ee
For example,  in the depth two we get
$$
    i^{|{\bf n}|-2} 
  \int_{(\R+i0)^2}  
\frac{e^{-ip_1w_1}}{{\mathfrak{sh}}^{a_1} (\pi p_1) \cdot {\mathfrak{sh}}^{b_1} (\pi \hbar p_1) }  
\frac{e^{-i(p_1+p_2)w_2}}{{\mathfrak{sh}}^{a_2} (\pi p_2) \cdot {\mathfrak{sh}}^{b_2}( \pi \hbar p_2)}
\frac{dp_1}{p_1^{n_1}}\frac{dp_2}{(p_1+p_2)^{n_2}}. 
$$
The functions ${\cal F}$ and $ {\cal I}$ are related by  
\be\la{ees}
 {{\cal F}}^\hbar_{{\bf a}, {\bf b}, {\bf n}}( \omega_1,  \omega_2, \ldots ,  \omega_{m}) = 
  {{\cal I}}^\hbar_{{\bf a}, {\bf b}, {\bf n}}( \omega_1-\omega_2,  \omega_2-\omega_3, \ldots ,  \omega_{m}).
\ee
Indeed, their integrands  differ only by the exponentials, related  by (\ref{popi1}).
Equivalently, we have 
\be\la{eessa}
 {{\cal F}}^\hbar_{{\bf a}, {\bf b}, {\bf n}}( \zeta_1, \ldots ,  \zeta_{m}) = 
  {{\cal I}}^\hbar_{{\bf a}, {\bf b}, {\bf n}}( \zeta_1+\ldots \zeta_m, \zeta_2+\ldots + \zeta_m,  \ldots ,  \zeta_{m}).
\ee

\subsubsection{Differential equations} 

\bp The quantum polylogarithms satisfy the differential equation
\be \la{iow}
d{{\cal F}}^\hbar_{{\bf a}, {\bf b}, {\bf n}}(\omega_1, ..., \omega_m) = \sum_{k=1}^m  {\cal F}^\hbar_{{\bf a} , {\bf b}, {\bf n}-{\rm 1}_k}(\omega_1, ..., \omega_m)d(\omega_k- \omega_{k+1}). 
\ee
\ep

\begin{proof} 
The introduced in Section \ref{secty} function
\be
\begin{split}
&  {\cal I}^\hbar_{{\bf a} , {\bf b}, {\bf n}}(\xi_1, ..., \xi_m) = \\ 
& i^{|{\bf n}|-m}\int_{(\R+i0)^m} {K}^\hbar_{{\bf a} , {\bf b}}(p_1, ..., p_m)
\frac{e^{-ip_1\xi_1} e^{-i(p_1+p_2)\xi_2} \ldots
e^{-i(p_1+...+p_m)\xi_m}dp_1 dp_2 \ldots dp_m}{p_1^{n_1}(p_1+p_2)^{n_2} \ldots 
(p_1+... +p_m)^{n_m}}.\\
 \end{split}
 \ee
evidently\footnote{since the variable  $\xi_k$ appears only in the exponential $e^{-ip_k\xi_k}$}  satisfies  the differential equations: 
\be \la{2.20.11.1}
d {\cal I}^\hbar_{{\bf a} , {\bf b}, {\bf c}}(\xi_1, ..., \xi_m) = 
\sum_{k=1}^m  {\cal I}^\hbar_{{\bf a} , {\bf b}, {\bf n}-{\rm 1}_k}(\xi_1, ..., \xi_m)d\xi_k.
\ee
Therefore the claim  follows from (\ref{ees}). \end{proof}

 \subsubsection{Complex conjugation} We claim that one has 
$$
\overline{{{\cal F}}_{{\bf a} , {\bf b}, {\bf n}}^\hbar(\omega_1, \ldots , \omega_m)}= (-1)^{|{\bf a}|+|{\bf b}|-m}{{{\cal F}}_{{\bf a} , {\bf b}, {\bf n}}^{\overline \hbar}(\overline \omega_1, \ldots , \overline \omega_m)}.
$$
Done by a change of 
variables $q=-\overline{p}$,  altering  the orientation of  $\R+i0$. Here is how it works in the depth one case.

\be
\begin{split}
& \overline{{\cal F}_{a,b,n}^\hbar(\omega)}= (-i)^{n-1}\int_{\overline{\R+i0}}
\frac{e^{i\overline p\overline{\omega}}}{{\mathfrak{sh}}^a(\pi \overline p){\mathfrak{sh}}^b(\pi \overline \hbar \overline p)}\frac{d\overline p}{\overline p^{n}} ~~ \stackrel{q=-\overline{p}}{=}\\
&-(-1)^{a+b} i^{n-1}\int_{{\R+i0}}
\frac{e^{-ip\overline{\omega}}}{{\mathfrak{sh}}^a(\pi p){\mathfrak{sh}}^b(\pi \overline  \hbar p)}\frac{dp}{p^{n}} = (-1)^{a+b+1} 
{\cal F}^{\overline \hbar}_{a,b,n}(\overline \omega).\\
\end{split}
\ee

\subsubsection{Shuffle relations} \la{shrel}
Quantum   polylogarithms satisfy shuffle relations, similar to the 
ones for the iterated integrals representing the 
  multiple polylogarithms.\footnote{It is  interesting that although 
quantum    polylogarithms do 
not seem to have an iterated integral presentation, 
they do satisfy the same shuffle relations as the iterated integrals.}
Namely, let us set
\be \la{abgf}
{\bf w}= (\omega_1, ..., \omega_m),   \quad {\bf u}= (u_1, ..., u_m), \quad {\bf t}=(t_1, ..., t_m). 
\ee
Consider the generating series 
in ${\bf u}$ whose coefficients  
are quantum polylogs with  indices  
${\bf n}$: 
\be \la{gfu}
\begin{split}
&{{\cal F}}^\hbar_{{\bf a} , {\bf b}}({\bf w}| {\bf u}):= 
\sum_{n_1, \ldots , n_m\geq 1}{{\cal F}}^\hbar_{{\bf a} , {\bf b}, {\bf n}}({\bf w}){u_1}^{{n}_1-1} \ldots {u_m}^{{n}_m-1}.\\
\end{split}
\ee
For example, 
$$
  {{\cal F}}^\hbar_{{a}, {b}}[w|u] =   i^{n-1} \int_{\R+i0} 
\frac{e^{-ip\omega}}{{\mathfrak{sh}}^{a}(\pi p){\mathfrak{sh}}^{b}(\pi \hbar p)}
~\frac{dp}{p-iu}. 
$$
Then make a   substitution $u_k:= t_1+ ... + t_k$ for $k=1, ..., m$:
$$
{{\cal F}}^\hbar_{{\bf a} , {\bf b}}[{\bf w}| {\bf t}^*]:= 
{{\cal F}}^\hbar_{{\bf a} , {\bf b}}({\bf w}| {t_1, t_1+t_2, \ldots , t_1+ \ldots + t_m}).
$$

\bl There is an integral presentation: 
$$
{{{\cal F}}}^\hbar_{{\bf a}; {\bf b}}[{\bf w}| {\bf t}^*]= i^{|{\bf n}|-m}\cdot 
 \int_{({\R+i0})^m}  
\prod_{k=1}^mK^\hbar_{a_k, b_k}(p_k; \omega_k)
\frac{dp_k}{(p_1+\ldots + p_k) - i(t_1+\ldots + t_k)}.
$$
\el

\begin{proof} Follows by applying the substitution $u_k:= t_1+ ... + t_k$ to the following identity:
\be \la{gens}
\sum_{n=1}^\infty \frac{dp_k}{(p_1+\ldots +p_k)^n}  (iu_k)^{n-1}=  \frac{dp_k}{(p_1+\ldots +p_k) - iu_k}. 
\ee  \end{proof}

Note that integrals (\ref{15f}) converge for any integers $n_i$, while the generating series use   $n_i \geq 1$. \vskip 2mm

Given   $
{\bf a}:= (a_1, ..., a_k)$ and $ {\bf a}':= (a_{k+1}, ..., a_{k+l})
 $ and a permutation 
$\sigma$ of the set $\{1, ..., k+l\}$, set
$
\sigma ({\bf a}{\bf a}'):= (a_{\sigma (1)},...,a_{\sigma (k+l) }).
$ 

\bt \label{shuffle} One has 
\begin{equation} \label{shuffle1}
{{{\cal F}}}^\hbar_{{\bf a}, {\bf b}}[{\bf w}| {\bf t}^*]\cdot {{{\cal F}}}^\hbar_{{\bf a'}, {\bf b'}}[{\bf w'}| {\bf t'}^*] = 
\sum_{\sigma \in \Sigma_{k,l}} {{{\cal F}}}^\hbar_{\sigma ({\bf a}{\bf a}'), \sigma ({\bf b}{\bf b}')}[\sigma ({\bf w}{\bf w}')
| \sigma ({\bf t}{\bf t}')^*].
\end{equation}
The sum is over the set  of 
all permutations shuffling $\{1, ..., k\}$ and $\{k+1, ..., k+l\}$.
\et

\begin{proof}  Follows immediately from the following identity  \cite[Lemma 2.12]{G01}:
\begin{equation} \label{A3}
\begin{split}
&\frac{1}{p_1( p_1 + p_2) ... ( p_1 + ... + p_k)} \cdot \frac{1}{p_{k+1}( p_{k+1}  + p_{k+2}) ... ( p_{k+1}   + ... + p_{k+l})} = \\
&\sum_{\sigma \in \Sigma_{k,l}} \frac{1}{p_{\sigma(1)}( p_{\sigma(1)} + p_{\sigma(2)} ) ... ( p_{\sigma(1)} + ... + p_{\sigma(k+l)} )}.\\
\end{split}
\end{equation} 
\end{proof}
For example, using  the identity 
$$
\frac{1}{p_1p_2} = \frac{1}{p_1(p_1+p_2)} + \frac{1}{p_2(p_1+p_2)}  
$$
 we have 
$$
{{\cal F}}^\hbar_{{a_1}, {b_1},{1}}(\omega_1)\cdot {{\cal F}}^\hbar_{{a_2}, {b_2}, {1}}(\omega_2) =
{{\cal F}}^\hbar_{(a_1, a_2), (b_1, b_2), (1,1)}(\omega_1, \omega_2) + {{\cal F}}^\hbar_{({a_2, a_1}), ({b_2, b_1}), (1,1)}(\omega_2, \omega_1). 
$$

\paragraph{\it Shuffle relations for the generating functions.} In addition to (\ref{abgf}), let us set
$$
  {\bf r}= (r_1, ..., r_m), \quad {\bf s}= (s_1, ..., s_m). 
$$
Generalizing (\ref{gfu}), we introduce the  quantum polylogarithm  generating series:
$$
{{\cal F}}^\hbar ({\bf w}| {\bf r}, {\bf s}, {\bf u}):= 
\sum_{a_i, b_i, n_i  >0}{{\cal F}}^\hbar_{{\bf a} , {\bf b}, {\bf n}}({\bf w})\prod_{k=1}^m r_k^{a_k-1}  {s_k}^{{b}_k-1}  {u_k}^{{n}_k-1}.
$$
One can rewrite this   using  the kernel  generating   series. 
\bl The kernel generating   function   is given by
$$
{\cal K}^\hbar(p; z| r,s):=   \sum_{a,b=1}^\infty   K^\hbar_{a,b}r^{a-1}s^{b-1} = \frac{e^{-ipz}}{({\mathfrak{sh}}(\pi p) - r)({\mathfrak{sh}}(\pi \hbar p)-s)}.
$$
\el

\bl The  quantum polylogarithm generating series   
 are given by the integrals  
 $$
 {{\cal F}}^\hbar ({\bf w}| {\bf r}, {\bf s}, {\bf u}):=  
  \int_{(\R+i0)^m} 
\prod_{k=1}^m{\cal K}^\hbar(p_k; \omega_k| r_k,s_k)
\frac{dp_k}{(p_1+\ldots + p_k)-iu_k}. 
$$
\el
 \begin{proof} Follows immediately using (\ref{gens}).
  \end{proof}
 
 For example,  the depth one quantum polylogarithm generating series are
 $$
 {{\cal F}}^\hbar (\omega| r, s, u):=  
 \int_{\R+i0} 
\frac{e^{-ip\omega}}{({\mathfrak{sh}}(\pi p) - r)({\mathfrak{sh}}(\pi \hbar p)-s)} 
\frac{dp}{ p -iu}. 
$$

Theorem \ref{shuffle}  immediately implies the following

\bt   One has 
$$
{{{\cal F}}}^\hbar[{\bf w}| {\bf r}, {\bf s}, {\bf t}^*]\cdot {{{\cal F}}}^\hbar [{\bf w'}| {\bf r}', {\bf s}', {\bf t'}^*] = 
\sum_{\sigma \in \Sigma_{k,l}} {{{\cal F}}}^\hbar[\sigma ({\bf w}{\bf w}')| {\sigma ({\bf r}{\bf r}'), \sigma ({\bf s}{\bf s}')}, 
 \sigma ({\bf t}{\bf t}')^*].
$$
\et

\subsubsection{Analytic continuation.} Quantum   polylogarithms 
have an analytic continuation to     
a cover 
of 
${\cal M}_{0, m+3}(\C)$ given by   
$$
(w_1, ..., w_m) \lra (\infty, -1, 0, e^{w_1}, ..., e^{w_m})\in {\cal M}_{0, m+3}(\C).
$$
The integral representation (\ref{2.27.11.2}) is convergent at the strip
$$
|{\rm Im}~\omega_i| < \pi a_i+\pi \hbar b_i. 
$$
We   use difference relations (\ref{2.20.11.2}) 
to extend it from that strip to $\C^m$, and argue by the 
induction on $|{\bf a}|+|{\bf b}|$. First, one checks formally 
 by induction that 
$$
\Delta^{(\omega_l)}_{i\pi \hbar}\Delta^{(\omega_k)}_{i\pi} 
{\cal I}^\hbar_{{\bf a}, {\bf b}, {\bf n}}({\bf w}) = 
\Delta^{(\omega_k)}_{i\pi}\Delta^{(\omega_l)}_{i\pi \hbar} 
{\cal I}^\hbar_{{\bf a}, {\bf b}, {\bf n}}({\bf w}).
$$
where each of the  sides is defined by applying twice 
  difference relations (\ref{2.20.11.2}).

  \subsubsection{An example: basic quantum   polylogarithms} \la{sec4.1}

We define the basic quantum   polylogarithms  by setting $a_j = b_j=1$. So in  the depth $m$   case
\be \la{153}
\begin{split}
&{{\cal F}}^\hbar_{n_1, ..., n_m}(\omega_1, ..., \omega_m):=  \\&
  i^{|{\bf n}|-m}\int_{(\R+i0)^m} 
\frac{e^{-ip_1\omega_1}}{{\mathfrak{sh}}(\pi p_1){\mathfrak{sh}}(\pi \hbar p_m)} \frac{dp_1}{p_1^{n_1} } \wedge \ldots \wedge \frac{e^{-ip_1\omega_1}}{{\mathfrak{sh}}(\pi p_m){\mathfrak{sh}}(\pi \hbar p_m)}\frac{dp_m}{(p_1+\ldots + p_m)^{n_m}}. \\
\end{split}
\ee

\bt \la{7.30.22} The function ${\cal F}_{n_1, ..., n_m}^\hbar(\omega_1, ..., n_m)$ enjoys the following properties: 

\begin{enumerate}

\item  Asymptotic expansion as $\hbar \to 0$:   
$$
 {\cal F}_{n_1, ..., n_m}^\hbar(\omega_1, ..., \omega_m) 
~\sim
 ~\frac{1}{(2\pi i \hbar)^m} \cdot   {\rm Li}_{{n_1+1}, \ldots , {n_m+1}}(e^{\omega_1-\omega_2}, e^{\omega_2-\omega_3},  \ldots, -e^{\omega_m}) 
 + \ldots 
$$


\item
Difference relations,  connecting  them with multiple polylogarithms:
\be \la{DIFRE33}
\begin{split}
&\Delta_{i\pi\hbar}^{(\omega_1)} \ldots \Delta_{i\pi\hbar}^{(\omega_m)}  {{\cal F}}_{n_1, ..., n_m}^\hbar(\omega_1, ..., \omega_m) ~=~ {\rm Li}_{n_1, ..., n_m}(e^{\omega_1-\omega_2}, e^{\omega_2-\omega_3}, \ldots , -e^{\omega_{m}}).  \\
&\Delta_{i\pi}^{(\omega_1)} \ldots \Delta_{i\pi}^{(\omega_m)}  {\cal F}_{n_1, ..., n_m}^\hbar(\omega_1, ..., \omega_m)  ~=~
\hbar^{|{\bf n}|-m} {\rm Li}_{n_1, ..., n_m}(e^{\frac{\omega_1-\omega_2}{\hbar}}, e^{\frac{\omega_2-\omega_3}{\hbar}}, \ldots , -e^{\frac{\omega_{m}}{\hbar}}).\\
\end{split}
\ee

 \item  
The value at $\hbar =1$: 
$$
{{\cal F}}^1_{n_1, ..., n_m}(\omega_1, ..., \omega_m) = {\bf P}_{{\bf 1}, {\bf 1}}\Bigl(\frac{\omega_1}{2\pi}. ..., \frac{\omega_m}{2\pi}\Bigr) 
\cdot {\rm Li}_{n_1, ..., n_m}\Bigl(e^{\omega_1}, ...,e^{\omega_m}\Bigr).
$$

\item  Distribution relations:
\be
\begin{split}
&r^{|{\bf n}|-m}{\cal F}_{n_1, ..., n_m}^{\frac{r}{s}\hbar}(r\omega_1, ..., r\omega_m) = \\
&\sum\limits_{\alpha_k=\frac{1-r}{2}}^{\frac{r-1}{2}} 
\sum\limits_{\beta_k=\frac{1-s}{2}}^{\frac{s-1}{2}} {\cal F}_{n_1, ..., n_m}^\hbar(\omega_1+\frac{2\pi i}{r}\alpha_1 
+ \frac{2\pi i\hbar }{s}\beta_1, \ldots , \omega_k+\frac{2\pi i}{r}\alpha_k 
+ \frac{2\pi i\hbar }{s}\beta_k ).\\
\end{split}
\ee

\item  The differential:
$$
d {\cal F}_{n_1, ..., n_m}^\hbar(\omega_1, ..., \omega_m) =   \sum_{j=1}^m
 {\cal F}_{n_1, ..., n_k-1, ..., n_m}^\hbar(\omega_1, ..., \omega_m) d(\omega_k-\omega_{k+1}). 
$$

\item  
$\hbar \longleftrightarrow 1/\hbar$ symmetry:
$$
{\cal F}_{n_1, ..., n_m}^\hbar(\omega_1, ..., \omega_m) = \hbar^{|{\bf n}|-m}{\cal F}_{n_1, ..., n_m}^{\frac{1}{\hbar}}\Bigl(\frac{\omega_1}{\hbar}, ..., \frac{\omega_m}{\hbar}\Bigr).
$$

\end{enumerate}
\et

 \begin{proof}   1)  
Using (\ref{36}) and (\ref{2.27.11.2a}), we get the leading term of the $\hbar \to 0$ asymptotic expansion:

\be \la{15a}
\begin{split}
& {{\cal F}}^\hbar_{n_1, ..., n_m}(\omega_1, ..., \omega_m)~\sim_{\hbar \to 0}~ \\
&\frac{i^{|{\bf n}|}}{(2\pi i \hbar)^m}  \int_{(\R+i0)^m} 
\frac{e^{-ip_1\omega_1}}{{\mathfrak{sh}}(\pi p_1)}\frac{dp_1}{p_1p_1^{n_1}} \wedge \ldots \wedge 
\frac{e^{-ip_m\omega_m}}{{\mathfrak{sh}}(\pi p_m)}
\frac{dp_m}{p_m (p_1+\ldots + p_m)^{n_m}} \\
& \stackrel{(\ref{2.27.11.2a})}{\sim}_{\hbar \to 0} \frac{1}{(2\pi i \hbar)^m} {\rm Li}_{{n_1}, \ldots , {n_m}}(e^{\omega_1-\omega_2},  e^{\omega_2-\omega_3},    \ldots,  -e^{\omega_m}) + \ldots .
\end{split}
\ee 
2)  For the second identity, set $q_i=p_i\hbar$, and 
use  integral  (\ref{2.27.11.2n}). The rest is straightforward. 
 \end{proof}

 \section{Quantum polylogarithms and  multiple $q-$polylogarithms}


Recall Definition \ref{MQP} of  the multiple $q-$polylogarithms: 
\be \la{MQPL}
\begin{split}
&{\rm Li}_{{\bf a}, {\bf n}}(x_1, ..., x_m; q):= \\
&\sum_{k_1,  \ldots ,  k_m>0}^{\infty}
\frac{x_1^{k_1}x_2^{k_2}\ldots x_m^{k_m}}{[k_1]_q^{a_1} [k_2]_q^{a_2}\ldots 
[k_m]_q^{a_m}\cdot k_1^{n_1} (k_1+k_2)^{n_2}\ldots (k_1+...+k_m)^{n_m}}.\\
\end{split}
\ee
If ${\bf a}=0$, we get the multiple polylogarithms 
$
{\rm Li}_{\bf n}(z_{1},\ldots , z_{m})   = {\rm Li}_{n_{1},\ldots ,n_{m}}(z_{1},\ldots ,z_{m}).
$

 Multiple $q-$polylogarithms  satisfy both the differential and difference equations:
 
  \paragraph{The 
differential.} 

Given an ${\bf x}=(x_1, \ldots , x_m)$, we set ${\bf x}^*:= (x_1...x_m, x_2...x_m, \ldots x_m)$. Then 
$$
d  {\rm Li}_{{\bf a}, {\bf n}}({\bf x}^*;q) =  \sum_{k=1}^m{\rm Li}_{ {\bf a}, {\bf n}- {\bf 1}_k}({\bf x}^*;q)d\log x_k. 
$$

\paragraph{The difference relation.} The $q-$difference operator defined by setting
 $$
 {\bf \Delta}_{x,q} f(x):=  f(qx) - f(q^{-1}x).
$$
We have the difference relations
\be \la{dife}
{\bf \Delta}_{x_k, q} {\rm L}_{{\bf a}, {\bf n}}({\bf x};q) =   {\rm L}_{ {\bf a}- {\bf 1}_k, {\bf n}}({\bf x};q).   
\ee

\subsection{Multiple $q-$polylogarithms by  the $q-$integration}

 \bd Given a  power series $f(x)$ and an integer $a\geq 0$,  the      $q-$integral  ${\Bbb I}_x^af(x)$ is:
$$
({{\Bbb I}}_x^af)(x):= (-1)^{{a}-1}\sum_{k\geq 0}{k+a-1\choose a-1}
 f(q^{2k+a}x).
 $$  
 \ed
 The name $q-$integral is justified by the following Lemma. 
 
 \bl \la{QDI} Let $a>0$. Then one has:
 $$
{\bf  \Delta}_{x,q} \circ {\Bbb I}_x^af(x)   = {\Bbb I}_x^{a-1}f (x).
 $$
  \el

 \begin{proof} Follows by a pretty standard calculation: 
\be
\begin{split}
&{\Bbb I}_x^af(qx) - {\Bbb I}_x^af(q^{-1}x) =\\
&(-1)^{{a}-1}\sum_{k\geq 0}{k+a-1\choose a-1}
\Bigl(f(q^{2k+a+1}x) - f(q^{2k+a-1}x)\Bigr) = \\
&(-1)^{{a}-1}\sum_{k\geq 0}
\Bigl({(k+1)+a-2\choose a-1}  
f(q^{2(k+1)+a-1}x) - {k+a-1\choose a-1}f(q^{2k+a-1}x)\Bigr) =\\
&(-1)^{{a}-1}\sum_{k> 0}
\Bigl({k+a-2\choose a-1} - {k+a-1\choose a-1}\Bigr) (q^{2k+a-1}x) - (-1)^{{a}-1}
f(q^{a-1}x)=\\
&(-1)^{{a}-2}\sum_{k\geq 0}
{k+a-2\choose a-2}  f(q^{2k+a-1}x) =\\
&  {\Bbb I}_x^{a-1}f(x).\\
\end{split}
\ee
\end{proof} 
 
 \bp \la{3.30.11.1} There is an equality of power series:
\be \la{3.18.11.1}
\begin{split}
&{\rm Li}_{ {\bf a}, {\bf n}}(x_1, ..., x_m; q)= ({\Bbb I}_{x_1}^{a_1}\ldots {\Bbb I}_{x_m}^{a_m}{\rm Li}_{{\bf n}})(x_1, \ldots,  x_m) = \\
&\sum_{k_1, ..., k_m\geq 0}(-1)^{|{\bf a}|-m}{k_1+a_1-1\choose a_1-1}\ldots {k_m+a_m-1\choose a_m-1}
{\rm Li}_{{\bf n}}(q^{2k_1+a_1}x_1, \ldots, q^{2k_m+a_m}x_m).\\
\end{split}
\ee
\ep

 \begin{proof} By Lemma \ref{QDI} and (\ref{dife}) we have
 $$
 {\bf \Delta}_{x_1, q}^{a_1}\circ \ldots \circ  {\bf \Delta}_{x_1, q}^{a_1}
 \Bigl({\rm Li}_{ {\bf a}, {\bf n}}(x_1, ..., x_m; q) -  ({\Bbb I}_{x_1}^{a_1}\ldots {\Bbb I}_{x_m}^{a_m}{\rm Li}_{{\bf n}})(x_1, \ldots,  x_m)\Bigr)=0.
  $$
   Both series vanish at $x_1=...=x_m=0$.
   \end{proof}

\paragraph{Examples}

\begin{enumerate}

\item {\it The $q$-polylogarithms} 
  are  power series in   $x$:
$$
{\rm Li}_{a,n}(x; q):= \sum_{k=1}^{\infty}\frac{x^k}{(q^k-q^{-k})^{a}~ k^n}, \qquad 
a,n \in \Z. 
$$
They satisfy both the differential and difference equations:
\be \la{DIFR44}
\begin{split}
&d {\rm Li}_{a, n}(x;q) =  {\rm Li}_{a, n-1}(x;q )d\log x. \\
&{\bf \Delta}{\rm Li}_{a, n}(qx;q)  =  {\rm Li}_{a-1, n}(x;q).\\
\end{split}
\ee
 
\item {\it Higher Pochhammer symbols $\Psi_{a+1}(x;q)$} are  the power series  given by the infinite products
$$
\Psi_{a+1}(x;q):=  
\prod_{n\geq 0}(1+q^{2n+1}x)^{(-1)^{a+1}{n+a\choose a}}.  
$$
For example, 
\be \la{dllog}
\begin{split}
&\Psi_1(x;q) = \frac{1}{(1+qx)(1+q^3x)(1+q^5x)(1+q^7x)\cdot  \ldots} , \\
&\Psi_2(x;q) = (1+qx)(1+q^3x)^2 (1+q^5x)^3(1+q^7x)^4\cdot \ldots .\\
\end{split}
\ee
They are the  unique power series in $x, q$  which satisfies the recursion 
$$
\frac{\Psi_a(qx;q)}{\Psi_a(q^{-1}x;q)} = \Psi_{a-1}(x;q), \qquad  \Psi_{0}(x;q) := 1+x. 
$$

\bp One has 
\be \la{PSI1}
\log \Psi_{a}(x;q) = -{\rm Li}_{a-1, 1}(-x; q)= -\sum_{k=1}^{\infty}\frac{(-x)^k}{(q^{k} - q^{-k})^{a-1} k}.
\ee
\ep 

\begin{proof}
Both power series satisfy the same difference equation, and equal 
to $0$ at $x=0$. 

Alternatively, here is a direct calculation for the classical case of $ \Psi_{1}(x;q)$. Formula (\ref{PSI1}) in this case is  the following identity: 
$$
{\rm Li}_{1, 1}(-x; q)=  \sum_{n\geq 0}\log(1+q^{2n+1}x).
$$
To prove it, we use the expansion   $\log(1+x) = -\sum_{k>0}\frac{(-x)^k}{k}$:
\be 
\begin{split}
&  \sum_{n\geq 0}\log(1+q^{2n+1}x) = - \sum_{n\geq 0} \sum_{k>0  } \frac{(-q^{2n+1}x)^k}{k}=\\
& -\sum_{k>0  } \sum_{n\geq 0}\frac{q^{2nk}q^k(-x)^k}{k} = - \sum_{k>0  } \frac{ q^k(-x)^k}{(1-q^{2k})\cdot k} = {\rm Li}_{1, 1}(-x; q).\\
\end{split}
\ee
\end{proof}

\item  Consider   a slight modification of the classical polylogarithm power series:
$$
{\rm L}_n(z) := -{\rm Li}_n(-z)= -\sum_{k>0}\frac{(-z)^k}{k^n}, ~~~|z|<1.
$$
  So  ${\rm L}_{1}(z) = \log(1+z)$ and  ${\rm L}_{0}(z) = \frac{z}{1+z}$. 
By Proposition \ref{3.30.11.1},  
$$
{\rm Li}_{a, n}(-x; q)= (-1)^{a}\sum_{k\geq 0}{k+a-1\choose a-1}{\rm L}_n(q^{2k+a}x).
$$

\item  It is interesting to compare  $q$-polylogarithms with the elliptic polylogarithms 
 \cite{BL}. The latter are obtained by the regularized weighted averaging  over $\Z$ of the 
classical polylogarithms, while the former are obtained by a similar 
weighted averaging but over the non-negative integers. For example, 
starting with $\log(1+x)$, the regularized averaging over $\Z$ delivers 
the logarithm of a theta function, while averaging over $\Z_{\geq 0}$ we get 
the negative of the logarithm of the $q$-exponential. 
\end{enumerate}

 \subsection{Quantum polylogarithms as  sums of $2^m$ companion $q-$polylogarithms} \la{sec3.2}

Recall that the quantum dilogarithm function can be written as a difference of two series:
\be \la{96}
  {{\cal F}}^\hbar_{1}(w) =  -{\rm Li}_{1, 1}(-e^{w}; q)  +{\rm Li}_{1,1}(-e^{w/\hbar}; q^\vee).
\ee
In Section \ref{sec3.2} we show that  quantum polylogarithms have similar presentation. 
We elaborate in detail  the case of    depth $m$  {\it basic quantum   polylogarithms}. 
We prove that they are sums of $2^m$ {companion    polylogarithm series}. One of them   is a quantum $q-$polylogarithm series,  another one 
is a  quantum $q^\vee-$polylogarithm series, and the other $2^{m}-2$  companion series are given by more general series. 
  
Recall the depth $m$   {\it basic quantum   polylogarithm}:   
$$
{{\cal F}}^\hbar_{n_1, ..., n_m}(\omega_1, ..., \omega_m):=  
i^{|{\bf n}|-m}  \int_{(\R+i0)^m} 
\prod_{k=1}^m\frac{e^{-i p_i\omega_k}}{{\mathfrak{sh}}(\pi p_k){\mathfrak{sh}}(\pi \hbar p_k)}
\frac{dp_k}{(p_1+\ldots +p_k)^{n_k}}. 
$$
Recall also the ${\cal I}-$variant (\ref{SECTSS}) of the basic quantum polylogarithms:
\be \la{154}
\begin{split}
 {{\cal I}}^\hbar_{n_1, ..., n_m}(w_1, ..., w_m) &:=  
 i^{|{\bf n}|-m} \int_{(\R+i0)^m} 
\prod_{k=1}^m \frac{e^{-i(p_1+...+p_k)w_k}}{{\mathfrak{sh}}(\pi p_k){\mathfrak{sh}}(\pi \hbar p_k)}\frac{dp_k}{(p_1+\ldots +p_k)^{n_k}}\\
&={{\cal F}}^\hbar_{n_1, ..., n_m}(w_1-w_2, w_2-w_3, ..., w_m).\\\end{split}
\ee

We calculate   integral (\ref{154})   as a sum over the residues. The sum splits   into a sum of $2^m$ series over the cones given by the direct sum of  $m$ copies 
 of the cones   
 $\Z_{>0}$ or $\hbar^{-1} \Z_{>0}$. 
 We call the resulting series {\it companion     polylogarithm series}.  
They are paramatrised by the sequences 
 $$
 \underline \varepsilon:= (\varepsilon_1, ..., \varepsilon_m), \quad \varepsilon_j \in \{1, \hbar^{-1}\}.
 $$
 For each such $ \underline \varepsilon$ we assign the {\it companion cone}:
$$
{\rm C}_{\underline \varepsilon}:=  \Z_{>0} \varepsilon_1 \oplus \ldots \oplus  \Z_{>0} \varepsilon_m \in \C^m.
$$
  Recall $q = e^{i\pi\hbar}$ and $q^\vee= e^{i\pi/\hbar}$. We also use a notation
  \be
  [k]_{q_{\varepsilon}} := \begin{cases} 
q^k-q^{-k}& \varepsilon=1\\
 (q^{\vee})^k-(q^{\vee})^{-k} & \varepsilon=\hbar^{-1}. 
\end{cases}.
 \ee
 
\bd The $\underline \varepsilon-$companion    polylogarithm series are given by
\be \la{CS}
\begin{split}
&{\rm Li}^\hbar_{{\bf a}, {\bf n}}({\underline \varepsilon}; w_1, ..., w_m) :=   \\
&\varepsilon_1 \ldots \varepsilon_m \sum_{k_1, ..., k_m>0}\frac{e^{k_1 \omega_1\varepsilon_1}\ldots e^{(k_1+...+k_m)\omega_m \varepsilon_m}}{[k_1]^{a_1}_{q_{\varepsilon_1}}  \ldots [k_m]^{a_m}_{ q_{\varepsilon_m}}  \cdot 
 (\varepsilon_1k_1)^{n_1}\ldots (\varepsilon_1 k_1+... +\varepsilon_m k_m)^{n_m}}.\\\end{split}
 \ee
 \ed
 The $\underline \varepsilon-$companion polylogarithm series (\ref{CS}) generalize   $q-$polylogarithm series (\ref{MQPL}). Indeed: 
 \be
 {\rm Li}^\hbar_{{\bf a}, {\bf n}}({\underline \varepsilon}; \omega_1, ..., \omega_m) =  \begin{cases} 
 {\rm Li}_{{\bf a}, {\bf n}}(e^{\omega_1}, ..., e^{\omega_m}; q)  & \mbox{if} ~\underline \varepsilon=(1, \ldots, 1),\\
  \hbar^{|{\bf n}|-m}{\rm Li}_{{\bf a}, {\bf n}}(e^{\omega_1/\hbar}, ..., e^{\omega_m/\hbar}; q^{\vee})   & \mbox{if} ~\underline \varepsilon=(\hbar^{-1}, \ldots, \hbar^{-1}). 
\end{cases}.
\ee
 
\bt \la{TH3.5} Assume $\hbar >0$. Assume that ${\rm Re}(\omega_i)<0$ and $|{\rm Im}(\omega_i)|<\pi$. Then we have 
\be
{\cal I}^\hbar_{{\bf n}}( {w_1}, \ldots,  {w_m}) = \sum_{\underline \varepsilon}{\rm Li}^\hbar_{{\bf 1}, {\bf n}}(\underline \varepsilon; w_1, w_2, \ldots , w_m).
\ee
\et

Before we proceed with the proof, let us elaborate  two examples.

\begin{enumerate}  \item  $m=1$. We get two companion cones:
$
\Z_{>0}$ and     $\hbar^{-1} \Z_{>0}. 
$ 
The related  companion series are  
  \be \la{CASS1}
\begin{split}
\Z_{>0}: \qquad &\sum_{k >0}\frac{(-1)^ke^{kw} }{[k]_q  ~ k^{n} } \ \ \ = \ \ {\rm Li}_{1,n}(-e^{w}; q).\\
 \hbar^{-1} \Z_{>0}: \qquad&\hbar^{n-1}\sum_{k>0}\frac{(-1)^ke^{kw/\hbar} }{[k]_{q^\vee} ~ k^{n}}= \hbar^{n-1}{\rm Li}_{1,n}(-e^{w/\hbar}; q^\vee).\\
\end{split}
 \ee
  So $ {{\cal I}}^\hbar_{  {n}}(w)$ is a sum of the two companion series (\ref{CASS1}):
  \be
  {{\cal I}}^\hbar_{  {n}}(w) =  {\rm Li}_{1,n}(-e^{w}; q)  + \hbar^{n-1}{\rm Li}_{1,n}(-e^{w/\hbar}; q^\vee).
  \ee
 When $n=1$ we recover formula (\ref{96}). 
 
\item  $m=2$.  We get four companion cones:
 $$
\Z_{>0} \oplus  \Z_{>0}, \quad \hbar^{-1}\Z_{>0} \oplus  \Z_{>0}, \quad \Z_{>0} \oplus \hbar^{-1} \Z_{>0}, \quad \hbar ^{-1}\Z_{>0} \oplus \hbar^{-1} \Z_{>0}. 
$$
The related companion series are:
  \be \la{CASS}
\begin{split}
\Z_{>0} \oplus  \Z_{>0}: \quad &\sum_{k_1, k_2>0}\frac{e^{-k_1w_1}e^{-k_2w_2}}{[k_1]_q [k_2]_q}\frac{1}{k_1^{n_1}(k_1+k_2)^{n_2}}.\\
\hbar^{-1}\Z_{>0} \oplus  \Z_{>0}: \quad &\hbar^{-1} \sum_{k_1, k_2>0}\frac{e^{-k_1w_1/\hbar}e^{-k_2w_2}}{[k_1]_{q^\vee} [k_2]_q }\frac{1}{( \hbar^{-1} k_1)^{n_1}(\hbar^{-1} k_1+  k_2)^{n_2}}.\\
\Z_{>0} \oplus \hbar^{-1} \Z_{>0}: \quad &\hbar^{-1} \sum_{k_1, k_2>0}\frac{e^{-k_1w_1}e^{-k_2w_2/\hbar}}{[k_1]_q  [k_2]_{q^{\vee}} }\frac{1}{k_1^{n_1}(k_1+\hbar^{-1} k_2)^{n_2}}.\\
 \hbar^{-1} \Z_{>0} \oplus \hbar^{-1} \Z_{>0}: \quad&\hbar^{-2}\sum_{k_1, k_2>0}\frac{e^{-k_1w_1/\hbar}e^{-k_2w_2/\hbar}}{[k_1]_{q^\vee} [k_2]_{q^{\vee}}}
 \frac{1}{(\hbar^{-1} k_1)^{n_1}(\hbar^{-1} k_1+\hbar^{-1} k_2)^{n_2}}.\\
\end{split}
 \ee
 So $ {{\cal I}}^\hbar_{  {n_1, n_2}}(w_1, w_2)$ is a sum of the four companion polylogarithm series (\ref{CASS}).
  \end{enumerate}\vskip 2mm
  
 \begin{proof} We elaborate   the case of the companion cone  $\hbar^{-1} \Z_{>0} \oplus \hbar^{-1} \Z_{>0}$. 
 Let ${\bf n}:=(n_1, n_2)$. Then  
$$
{{\cal I}}^\hbar_{  {\bf n}}(w_1, w_2):=  
 \int_{({\R+i0})^2}  
\frac{e^{-ip_1w_1}}{{\mathfrak{sh}} (\pi p_1){\mathfrak{sh}} (\pi \hbar p_1)}
\frac{e^{-i(p_1+p_2)w_2}}{{\mathfrak{sh}} (\pi p_2){\mathfrak{sh}} (\pi \hbar p_2)}
\frac{dp_1}{p_1^{n_1}}\frac{dp_2}{(p_1+p_2)^{n_2}}. 
$$
The  contribution of the residues at   $ p_1=  i k_1 /\hbar,   p_2=  i k_2 /\hbar $  where $k_1, k_2>0$   gives 
\be
\begin{split}
&\hbar^{-2}\sum_{k_1, k_2>0}  
\frac{ (e^{w_1 /\hbar})^{k_1}}{[k_1]_{q^\vee} }
\frac{ (-e^{w_2/\hbar})^{k_1+k_2}}{[k_2]_{q^\vee} }
\frac{1}{(\hbar^{-1} k_1)^{n_1} ~(\hbar^{-1}(k_1+k_2))^{n_2}} \\&= 
 \hbar^{|{\bf n}|-2}{\rm Li}_{{\bf 1},  {\bf n}}(e^{w_1/\hbar}, -e^{w_2/\hbar}; q^\vee). \\
\end{split}
\ee
\end{proof}

The analog of Theorem \ref{TH3.5} for arbitrary quantum polylogarithms is obtained  by a similar residue calculation. Since 
 the function 
$
\frac{1}{{\mathfrak{sh}}^{a_s} (\pi p_s){\mathfrak{sh}}^{b_s} (\pi  \hbar p_s) } 
$
   has zeros of higher order at $p_s= ik_s, i  k_s/\hbar$ where $k_s>0$, we get   sums of companion series  multiplied by powers of $\omega_s$. 
   Calculating the residues at $p_s= i  k_s/\hbar$ we encounter the following derivatives, evaluated then at $p_s= i  k_s/\hbar$:
   \be
\begin{split}
&
   \left(\frac{d}{dp_s}\right)^{b_s-1}  \frac{1}{{\mathfrak{sh}}^{a_s} (\pi   p_s)  } 
   \prod^m_{j=s}\frac{e^{-i(p_1+...+p_j)z_j}}{  (p_1+...+p_j)^{n_j}}.\\
     \end{split}
     \ee
For the residues at $p_s= i k_s$ we get  similar derivatives, with $a_s$ switched with $b_s$,   at $p_s= i k_s$:

 \section{Depth one examples}

 We consider integrals the depth one integrals for different countours 
 $\alpha$:
$$
 i^{n-1}\cdot  \int_{\alpha} \frac{e^{-ipz}}{ {\mathfrak{sh}}^a (\pi p){\mathfrak{sh}}^b (\pi \hbar p) }\frac{dp}{p^{n}}.
$$
 If  $\alpha:= \alpha_0$  is a small  {\it counterclockwise} oriented loop around zero,   we get   polynomials in $z$,   generalizing Bernoulli polynomials. 
If  $\alpha:=\R+i0$, we get  the depth one  quantum  polylogarithms.

\subsection{Quantum Bernoulli polynomials} \la{sec4.2}

Recall the Bernoulli polynomials $B_n(x)$: 
$$
\frac{te^{tx}}{e^t-1} = \sum_{n=0}^\infty B_n(x)\frac{t^n}{n!}.
$$

\bd Quantum Bernoulli polynomials are polynomials in $\omega$ and $\hbar^{\pm 1}$ given by
$$
B^\hbar_{a,b, n}(\omega):= i^{n-1}\int_{\alpha_0} \frac{e^{-ip\omega}}{ {\mathfrak{sh}}^a (\pi p){\mathfrak{sh}}^b (\pi \hbar p) }\frac{dp}{p^{n}}.
$$
\ed

To state the properties of the polynomials   $B^\hbar_{a,b,n}(\omega)$ we need   
  a polynomial
$$
Q_{m}(\omega) = \frac{\Bigl(\omega - \pi i (m-1)\Bigr)\Bigl(\omega-\pi i (m-3)\Bigr) \cdot \ldots \cdot \Bigl(\omega- \pi i (1-m)\Bigr)}{ (2\pi i)^mm!}.
$$
It is the unique degree $m$ polynomial  with the following two properties: 

\begin{itemize}

\item It satisfies   difference relations
\be \la{2.2.11.1}
\Delta_{i\pi} Q_m(\omega) = Q_{m-1}(\omega), \quad Q_0(\omega) =1.
\ee
 
\item The roots of $Q_m(\omega)$ form an arithmetic progression with the step $2\pi i$, centered at $0$. 

 \end{itemize}
 
  The weight of $Q_{m}(\omega)$ is $0$.  For example
 \be \la{Q1}
 Q_0(\omega) =1, ~~~~Q_1(\omega) =\frac{\omega}{2\pi i}, ~~~~Q_2(\omega) =\frac{(\omega - i\pi)(\omega+i\pi)}{2!\cdot (2\pi i)^2}.
\ee

\bt The quantum Bernoulli polynomials  $B^{\hbar}_{a,b,n}(\omega)$
have  the following properties.
 
\begin{enumerate}
\item $B^\hbar_{a,b, n}(\omega)$ is a  polynomial   
in $\omega$ of the degree $a+b+n-1$.  

\item  Differential and difference equations:
\be
\begin{split}
d B^{\hbar}_{a,b, n}(\omega ) &=  B^{\hbar}_{a,b, n-1}(z )d\omega. \\
\Delta_{i\pi\hbar}B^{\hbar}_{a,b, n}(\omega)   &=
B^{\hbar}_{a,b-1, n}(\omega), \\
\Delta_{i\pi}B^{\hbar}_{a,b, n}(\omega) &=
B^{\hbar}_{a-1,b, n}(\omega).  \\
\end{split}
\ee

\item  Asymptotic expansion when $\hbar \rightarrow 0$:
$$ B^\hbar_{a,b, n}(\omega) \   \sim_{\hbar \to 0} \   \frac{1}{(2\pi i\hbar)^b}  B_{a, 0, b+n}(\omega) + \ldots .
$$
 
\item The value at $\hbar =1$:
$$B_{a,b,n}^1(\omega) = B_{a+b, 0, n}(\omega).
$$

\item  Relation with  Bernoulli polynomials $B_n(\omega)$ and  polynomials $Q_n(\omega)$:
\be \la{2.2.11.1a}
\begin{split}
&B^{\hbar}_{1,0,n}(\omega ) =  \frac{(2\pi i)^{n-1}}{(n-1)!}B_{n-1}\left(\frac{\omega}{2\pi i}+ \frac{1}{2}\right),\\
&B^{\hbar}_{a, 0, 1}(\omega) = Q_{a-1}(\omega).\\
\end{split}
\ee

\item Modular property, or   $\hbar \leftrightarrow 1/\hbar $ symmetry: for any $\hbar \in \C^\times$ one has:
$$
B^\hbar_{a,b,n}(\omega) = \hbar^{n-1}B^{1/\hbar}_{b,a, n}(\omega/\hbar).
$$

\item Complex conjugation:
$$
\overline{B^\hbar_{a,b, n}(\omega)}= (-1)^{a+b+1}{B^{\overline \hbar}_{a,b, n}(\overline \omega)}.
$$

\item The $\omega \longleftrightarrow -\omega$ symmetry:
$$
B^\hbar_{a,b, n}(\omega) = (-1)^{a+b+n}B^{\hbar}_{a,b,n}(-\omega ). 
$$

\item The generating function 
$$
B^\hbar(\omega|r,s,u):=  \int_{\alpha_0} \frac{e^{-ip\omega}}{({\mathfrak{sh}}(\pi p) - r)({\mathfrak{sh}}(\pi \hbar p)-s)}\frac{dp}{p -i u}.
$$

\end{enumerate}
\et

\begin{proof}  We present an argument only if it is not totally straightforward.  

1) A residue calculation. 

5) The first claim is an easy calculation. The second is proved in the following Lemma. 

\bl \la{2.13.11.1} One has for $m \geq 0$\footnote{The $-i$ factor is just  the  factor $i^{n-1}$ in the case  $n=0$.}
\be \la{2.12.11.1a}
-i\int_{\alpha_0}\frac{e^{-ip\omega}}{{\mathfrak{sh}}^{m}(\pi p)}dp= 
Q_{m-1}(\omega).
\ee
\el

For example, 
$$
-i\int_{\alpha_0}\frac{e^{-ip\omega}}{{\mathfrak{sh}}(\pi p)}dp = 1, \qquad
-i\int_{\alpha_0}\frac{e^{-ip\omega}}{{\mathfrak{sh}}^2(\pi p)}dp = \frac{\omega}{2\pi i}.
$$
 
\begin{proof}  Integral (\ref{2.12.11.1a}) satisfies recursion (\ref{2.2.11.1}), which determines 
each next one uniquely up to a constant.  
Furthermore, $Q_m(\omega) = (-1)^{m}Q_m(-\omega)$, which tells that $Q_{2k+1}(\omega)$ is divisible by $\omega$.  
This, however, does not complete the proof, so we give a proof based on a different idea. 
Set 
$$
{\rm I}_m(\omega):= -i\int_{\alpha_0}\frac{e^{-ip\omega}}{{\mathfrak{sh}}^{m}(\pi p)}dp. 
$$
Then one has a recursion
\be \la{rec11}
{\rm I}_{m+1}(\omega)=  \frac{\omega-i\pi (m-1)}{2\pi i m}{\rm I}_{m}(\omega+i\pi ).
\ee
Indeed, integrating by parts we get
$$
 \int_{\alpha_0} e^{-ip\omega}\frac{d}{dp}{\mathfrak{sh}}^{-m}(\pi p)dp = \omega{\rm I}_m(z).
$$
Since
$$
-\frac{d}{dp}\mathfrak{sh}^{-m}(\pi p) = \frac{\pi m\cdot (e^{\pi p}+e^{-\pi p})}{\mathfrak{sh}^{m+1}(\pi p)} =
\frac{\pi m}{\mathfrak{sh}^{m}(\pi p)} + \frac{2\pi m \cdot e^{-\pi p}}{\mathfrak{sh}^{m+1}(\pi p)}, 
$$
we get
$
 (\omega-i\pi m)\cdot {\rm I}_m(\omega) = 2\pi i m \cdot {\rm I}_{m+1}(\omega-i\pi). 
$ 
 This is equivalent to (\ref{rec11}). Therefore 
$$
{\rm I}_{m+1}(\omega) = Q_m(z){\rm I}_{1}(\omega+i\pi m )= Q_m(\omega).
$$
\end{proof}

 6) Done by a change of 
variables $q=p/\hbar$, preserving the isotopy class of the contour $\alpha_0$. \vskip 1mm

7) Done by a change of 
variables $q=-\overline{p}$,  altering  the orientation of  $\alpha_0$:

\be
\begin{split}
& \overline{{B}_{a,b,n}^\hbar(\omega)}= (-i)^{n-1}\int_{\overline{\alpha_0}}
\frac{e^{i\overline p\overline{\omega}}}{{\mathfrak{sh}}^a(\pi \overline p){\mathfrak{sh}}^b(\pi \overline \hbar \overline p)}\frac{d\overline p}{\overline p^{n}} ~~ \stackrel{q=-\overline{p}}{=}\\
&(-1)^{a+b+1} i^{n-1}\int_{{\alpha_0}}
\frac{e^{-ip\overline{\omega}}}{{\mathfrak{sh}}^a(\pi p){\mathfrak{sh}}^b(\pi \overline  \hbar p)}\frac{dp}{p^{n}} = (-1)^{a+b+1} 
{B}^{\overline \hbar}_{a,b,n}(\overline \omega).\\
\end{split}
\ee
The extra $-$ sign amounts to the change of the contour orientation. 

8) It is obtained by a change of 
variables $q=-p$. It does not change the contour $\alpha_0$. 
\end{proof}

\subsection{Depth one quantum   polylogarithms}  \la{Sec2}

They  are   the following integrals:
$$
{{\cal F}}^\hbar_{a,b,n}(\omega) := i^{n-1}\int_{\R+i0}\frac{e^{-ip\omega}}{ {\mathfrak{sh}} (\pi p)^a
{\mathfrak{sh}} (\pi \hbar p)^b } \frac{dp}{p^{n}}, ~~~~a,b\in \Z_{\geq 0}, \ n \in \Z.
$$

The next Lemma tells that they reduce  to the classical polylogarithms when $b=0$.
\bl \la{QMP2} i) One has for $a > 0$, $m\geq 0$:
\be \la{2.16.11.1x}
i^{-m-1}\int_{\R+i0}\frac{e^{-ip\omega}}{{\mathfrak{sh}}^{a}(\pi p)} p^{m}dp = \left(\frac{d}{d\omega}\right)^{m}\Bigl(Q_{a-1}(\omega) 
{\rm Li}_{0}(e^{\omega+\pi i  a })\Bigr).
\ee
In particular,   it is a single-valued meromorphic function in $z$

ii) One has for $a > 0$, $n>  0$:
\be \la{2.16.11.1}
\begin{split}
&i^{n-1}\int_{\R+i0}\frac{e^{-ip\omega}}{{\mathfrak{sh}}^{a}(\pi p)}\frac{dp}{p^{n}} \ = \ i^{n}\sum_{k\geq 0} {n+k\choose k} \left(-\frac{d}{d\omega}\right)^{k}Q_{a-1}(\omega)\cdot 
{\rm Li}_{n+k}\left(e^{\omega+\pi i  a}\right).\\
\end{split}
\ee
\el

Formulas (\ref{2.16.11.1}) look simpler for the generating series ${\rm Li}(x; t):= \sum_{n>0}{\rm Li}_n(x)t^{n-1}$
 \be \la{2.16.11.1aa}
\begin{split}
&\sum_{n>0} t^{n-1}\int_{\R+i0}\frac{e^{-ip\omega}}{{\mathfrak{sh}}^{a}(\pi p)}\frac{dp}{p^{n}} \cdot t^{n-1}\ = \ i^{n} \left(1 + t^{-1} \frac{d}{dw}\right)^{-1}Q_{a-1}(\omega)\cdot 
{\rm Li}\left(e^{\omega+\pi i  a}; t\right).\\
\end{split}
\ee 

\paragraph{ Examples.} 1. Note that ${\rm Li}_{0} (e^{\omega} )= \frac{d}{d\omega}{\rm Li}_1(e^\omega) =\frac{e^{\omega}}{1-e^{\omega}}$. Then one has 
\be \la{2.12.11.1}
\begin{split}
&i^{-1}\int_{\R+i0}\frac{e^{-ip\omega}}{{\mathfrak{sh}}^{a}(\pi p)}dp ~= ~Q_{a-1}(\omega){\rm L}_0(e^{\omega+\pi i a}) ~=~
 Q_{a-1}(\omega)\frac{e^{\omega+\pi i a}}{1- e^{\omega+\pi i a}}.\\
\end{split}
\ee
For example, 
\be \la{2.12.11.1*}
\begin{split}
& i^{-1}\int_{\R+i0}\frac{e^{-ip\omega}}{{\mathfrak{sh}}(\pi p)}dp = \frac{-e^{\omega}}{1+e^{\omega}}, \\
&i^{-1}\int_{\R+i0}\frac{e^{-ip\omega}}{{\mathfrak{sh}}^2(\pi p)}dp = \frac{\omega}{2\pi i}\frac{e^{\omega}}{1-e^{\omega}},\\ 
&i^{-1}\int_{\R+i0}\frac{e^{-ip\omega}}{{\mathfrak{sh}}^3(\pi p)}dp = \frac{\omega^2+\pi^2}{2!(2\pi i)^2}\frac{-e^{\omega}}{1+e^{\omega}}.\\ 
\end{split}
\ee 

2. Formula (\ref{2.16.11.1x}) for $m=0$ is 
\be \la{2.16.11.1d}
i^{-1}\int_{\R+i0}\frac{e^{-ip\omega}}{{\mathfrak{sh}}^{a}(\pi p)}\frac{dp}{p} = \sum_{k\geq 0}\left(-\frac{d}{d\omega}\right)^{k}Q_{a-1}(\omega) \cdot {\rm Li}_{k}\left(e^{\omega+\pi i  a}\right).
\ee

\begin{proof} i) We start with the  $m=0$ case. Let us calculate  integral (\ref{2.12.11.1}) using the residue theorem, 
assuming that ${\rm Re}(\omega)<0$, using 
the  rectangular contour ${\Omega}_{N}$  a bit over the real axis. 
The residues are at the points $ik$, $k>0$. The residue at $p=ik$ equals to the residue at $p=0$ multiplied by $(-1)^{ka}e^{k\omega}$.  Lemma \ref{2.13.11.1} calculaties the integral around $p=0$. So we get 
$$
Q_{a-1}(\omega) \cdot \sum_{k>0}(-1)^{ka}e^{k\omega} = Q_{a-1}(\omega) \cdot \sum_{k>0}e^{k(\omega+i\pi a)} = 
Q_{a-1}(\omega)\frac{e^{(\omega+i\pi a)}}{1-e^{(\omega+i\pi a)}}.
$$
Formula (\ref{2.16.11.1x}) follows from this by differentiating by $\omega$.

ii) We prove   (\ref{2.16.11.1aa}) applying $1+t^{-1}\frac{d}{d\omega}$ to the left hand side, using the fact that 
$$
 i^n \frac{d}{d\omega} \int_{\R+i0}\frac{e^{-ip\omega}}{{\mathfrak{sh}}^{a}(\pi p)}\frac{dp}{p^{n}} =  i^{n-1}\int_{\R+i0}\frac{e^{-ip\omega}}{{\mathfrak{sh}}^{a}(\pi p)}\frac{dp}{p^{n-1}} .
 $$

\end{proof}

\bt \label{qdilog} The depth one quantum polylogarithm ${{\cal F}}^\hbar_{a,b,n}(\omega)$
has the following features:
\vskip 2mm
\begin{enumerate}

\item  Differential and difference equations:
\be \la{PRO11}
\begin{split}
d {{\cal F}}_{a,b, n}^\hbar(\omega) &=  {{\cal F}}_{a,b, n-1}^\hbar(\omega) d\omega,\\
\Delta_{i\pi  \hbar}{{\cal F}}^\hbar_{a,b, n}(\omega)  &=
{{\cal F}}_{a,b-1,c}^\hbar(\omega), \\
\Delta_{i\pi}{{\cal F}}^\hbar_{a,b, n}(\omega) &=
{{\cal F}}_{a-1,b,n}^\hbar(\omega). \\ 
\end{split}
\ee

\item The limit  when $\Re z \rightarrow -\infty$,  taken along a line parallel to the real
axis: 
$$
\lim_{\Re z \rightarrow -\infty}{\cal F}_{a,b,n}^\hbar(z)=0. 
$$

\item   Let $a, b\geq 0$ and $n \geq 1$. Then 
${{\cal F}}^\hbar_{a,b,n}(\omega)$ is a  multivalued analytic function  with the singularities 
at the two integral positive cones:
$$
\left \{\pm \pi i \Bigl((2m+a) + 
(2n+b)\hbar)\Bigr)~|~m,n \in {\mathbb \Z_{\geq 0}}\right \}. 
$$

\item  Asymptotic expansion when $\hbar \rightarrow 0$:
\be
\begin{split}
& {{\cal F}}^\hbar_{a,b,n}(\omega) ~\sim_{\hbar \rightarrow 0}  ~
  \frac{i^{n}}{(2\pi \hbar)^b} \sum_{k\geq 0} {n+k\choose k} \left(-\frac{d}{d\omega}\right)^{k}Q_{a-1}(\omega)\cdot 
{\rm Li}_{b+n+k}\left(e^{\omega+\pi i  a}\right) + \ldots.\\
\end{split}
\ee

\item The value at $\hbar=1$:
\be
\begin{split}
& {{\cal F}}_{a,b,n}^1(\omega) = 
  i^{n+b}\sum_{k\geq 0} {n+k\choose k} \left(-\frac{d}{d\omega}\right)^{k}Q_{a+b-1}(\omega) {\rm L}_{n+k+b}\left(e^{\omega+\pi i (a+b)}\right).\\
\end{split}
\ee

\item Complex conjugation:
$$
\overline{{{\cal F}}_{a,b,n}^\hbar(\omega)}= (-1)^{a+b-1}{{{\cal F}}_{a,b,n}^\hbar(\overline \omega)}.
$$

\item The $\omega \longleftrightarrow -\omega$ symmetry:
$$
{{\cal F}}_{a,b,n}^\hbar(\omega) + (-1)^{a+b+n-1}{{\cal F}}_{a,b,n}^\hbar(-\omega)= B_{a,b, n}(\omega; \hbar). 
$$

\item The $\hbar \longleftrightarrow 1/\hbar $ symmetry:
\be\la{MODP11}
{{\cal F}}^\hbar_{a,b,n}(\omega) =  \hbar^{n-1} {{\cal F}}_{b,a,n}^{1/\hbar}(\omega\hbar).
\ee

\item Distribution relations: 
$$ 
 r^{n-1}{\cal F}^{\frac{r}{s}\hbar}(r\omega) = \prod\limits_{\alpha=\frac{1-r}{2}}^{\frac{r-1}{2}}
\prod\limits_{\beta=\frac{1-s}{2}}^{\frac{s-1}{2}} {\cal F}^\hbar(\omega+\frac{2\pi i}{r}\alpha
+ \frac{2\pi i\hbar }{s}\beta ). 
$$

\end{enumerate}

\et

\begin{proof}    
  We provide the arguments only when they are not evident. 
  
    3) If $\hbar \in \R_{>0}$  the integral converges when   
  $|{\rm Im} z| < 1+\hbar$. The claim follows from   recursions  (\ref{PRO11}), and uses Property 2) for the normalization. 
  Precisely, let $n=1$. Then by Lemma \ref{QMP2}
  
1) The function 
  ${\cal F}^\hbar_{a,0,1}(z)$  has simple poles at the rays $\pm \pi i (2\Z_{\geq 0}+a)$.  
  
2) The function  ${\cal F}^\hbar_{0,b,1}(z)$  
  has simple poles at the rays   $\pm \pi i \hbar (2\Z_{\geq 0}+b)$.

  Note that the roots of polynomials $Q_{a-1}(z)$ and $Q_{b-1}(z/\hbar)$ kill the poles at the    centered at $0$ segments of lengths $a-2$ and  respectively $b-2$, 
  with the steps $2\pi i$ and $2\pi i \hbar$. 
 
  The case $n>1$ is obtained by the integration in $z$, and thus follows from the $n=1$ case. 
  The   case when $ \hbar \in \C - (-\infty, 0]$ follows by an analytic continuation. 
  
    4) Follows  by (\ref{2.16.11.1}) from
 \be
\begin{split}
& {{\cal F}}^\hbar_{a,b,n}(\omega) ~\sim_{\hbar \to 0}~ \frac{i^{n-1}}{(2\pi 
  \hbar)^b}\int_{\R+i0}\frac{e^{-ip\omega}}{ {\mathfrak{sh}} (\pi p)^a} \frac{dp}{p^{b+n}} + \ldots \\
\end{split}
\ee

  5) Follows from (\ref{2.16.11.1}). 
 
 7) Change the  
variables $q=- {p}$. It   changes the  integration contour $\gamma$  to $-\gamma$. Their sum  is a clockwise contour  around $0$. 
\end{proof}

\paragraph{Example: Basic depth one quantum polylogarithms.} They are given by the integrals 
$$
{{\cal F}}^\hbar_{n}(\omega) := i^{n-1}\int_{\R+i0}\frac{e^{-ip\omega}}{ {\mathfrak{sh}} (\pi p) 
{\mathfrak{sh}} (\pi \hbar p)  } \frac{dp}{p^{n}}, ~~~~  n \in \Z, 
$$ 
and have  the following properties: 

\begin{enumerate}

\item  Asymptotic relation to the $n$-logarithm 
$$
 {\cal F}_n^\hbar(\omega) ~\sim_{\hbar \rightarrow 0} ~\frac{{\rm L}_{n+1}(-e^{\omega})}{2\pi 
  \hbar}. 
$$
\item  The differential and difference relations:
\be \la{DIFRE33}
\begin{split}
&d {\cal F}_{n}^\hbar(\omega) =  {\cal F}_{n-1}^\hbar(\omega)d\omega \\
&\Delta_{i\pi \hbar}{\cal F}_{n}^\hbar(\omega) =
{\rm L}_{n}(-e^\omega), \\
&\Delta_{i\pi}{\cal F}_{n}^\hbar(\omega)  =
\hbar^{n-1}{\rm L}_{n}(-e^{\omega/\hbar}).  \\
\end{split}
\ee

 \item  
Modular property:
$$
{\cal F}_{n}^\hbar(\omega) + \hbar^{1-n}{\cal F}_{n}^{-1/\hbar}(\omega/\hbar) = 0.
$$
 Equivalently, the generating series
$
{\cal F}(\omega; \hbar; t) := \sum_{m=1}^\infty{\cal F}_{m}^\hbar(\omega)t^{m-1}
$
satisfy
$$
{\cal F}(\omega; \hbar; t)   + {\cal F}(\omega/\hbar; -1/\hbar;  t/\hbar)=0.
$$
 
\item  
Relation with  $q$-polylogarithms:
\be \la{PROP22}
{\cal F}_{n}^\hbar(\omega) ~=~
{\rm L}_{1, n}(e^{-\omega}; e^{i\pi \hbar})  - \hbar^{n-1}{\rm L}_{1, n}(e^{-\omega/\hbar}; {e^{i\pi /\hbar}}). 
\ee
Equivalently, the generating series ${\rm L}(e^z; e^{i\pi \hbar}; t) := \sum_{m=1}^\infty{\rm L}_{1, m}(e^z; e^{i\pi \hbar}) t^{m-1}$ satisfy
$$
{\cal F}(\omega; \hbar; t)  = {\rm L}(e^\omega; e^{i\pi \hbar}; t)  - {\rm L}(e^{\omega/\hbar}; e^{-i\pi /\hbar}; \hbar t).
$$

 \item  Distribution relations:
$$ \sum\limits_{\alpha=\frac{1-r}{2}}^{\frac{r-1}{2}} 
\sum\limits_{\beta=\frac{1-s}{2}}^{\frac{s-1}{2}} {\cal F}_n^\hbar(\omega+\frac{2\pi i}{r}\alpha 
+ \frac{2\pi i\hbar }{s}\beta ) = 
r^{n-2}{\cal F}_n^{\frac{r}{s}\hbar}(r\omega).
$$ 
\end{enumerate}

 \section{Concluding remarks}
 
 Scattering amplitudes in the ${\cal N}=4$ SUYM theory can be expressed via polylogarithms and their generalizations. For example, the MHV $n$ particles ${\rm L}-$loop 
 scattering amplitudes are weight $2{\rm L}$ functions on the configuration space ${\rm Conf}_n(\C{\Bbb P}^3)$ of collections of 
 $n$ points in $\C{\Bbb P}^3$, considered modulo the diagonal action of the group ${\rm PGL}_4$. It is  invariant under the cyclic shift of the points. 
 The space ${\rm Conf}_n({\Bbb P}^3)$ carries canonical cluster Poisson structure, invariant under the cyclic shift, 
  and the related space  ${\rm Conf}_n(\C^4)$ carries a cluster $K_2-$variety structure. 
 \vskip 1mm
 What is  role the cluster Poisson structure on ${\rm Conf}_n({\Bbb P}^3)$  for  the  scattering amplitudes? 
 \vskip 1mm
  
 Any cluster Poisson variety $\mathscr{X}$ admits cluster quantization \cite{FG3}, where the quantised algebra of functions ${\cal O}_q(\mathscr{X})$  acts 
 by unbounded operators in a cluster Hilbert space ${\cal H}_{\mathscr{X}}$. The Hilbert space ${\cal H}_{\mathscr{X}}\otimes \overline {\cal H}_{\mathscr{X}}$ for a given cluster coordinate system ${\bf c}$ is realised as the  Hilbert space 
 ${\rm L}_2(\mathscr{A}(\R_{>0}))$ of functions on the space of real positive points of the cluster $K_2-$variety $\mathscr{A}(\R_{>0})$. These Hilbert spaces for different cluster coordinate systems are related by the quantum dilogarithm intertwiners. 
 
 Suppose that the asymptotic expansion as $\hbar \to 0$ of  a vector $\varphi^{\bf c}_\hbar$ in the Hilbert space for a single cluster coordinate system ${\bf c}$ is written as 
 \be \la{1/1/2026}
 \varphi^{\bf c}_\hbar \sim e^{F_{\bf c}(a, \hbar)/\hbar}, \ \ \ \mbox{where 
 $F_{\bf c}(a, \hbar)$ is a function on 
 $\mathscr{A}(\R_{>0})$ depending on $\hbar$}.
 \ee
   The vectors $\varphi^{\bf c}_\hbar$ and $\varphi^{\bf c'}_\hbar$ for two different clusters ${\bf c}$ and ${\bf c}'$ are related by the quantum 
 dilogarithm intertwiners. Therefore the stationary phase method shows that the $\hbar \to 0$ asymptotics of the 
 vectors $\varphi^{\bf c}_\hbar$ in any cluster coordinate system ${\bf c}$ can be written in the form (\ref{1/1/2026}). 
 Moreover the functions   $F_{\bf c}(a, \hbar)_{|\hbar=0}$ and  $F_{\bf c'}(a, \hbar)_{|\hbar=0}$ for any 
 two clusters ${\bf c}$ and ${\bf c}'$ are the same functions on $\mathscr{A}(\R_{>0})$, expressed in the cluster coordinates for the clusters ${\bf c}$ and ${\bf c}'$. \vskip 1mm
 
 I suggest that  the scattering amplitudes should have an $\hbar-$deformation, 
  becoming vectors ${\cal A}^\hbar_{n, {\rm L}}$ in the space of cluster distributions. These vectors  should be expressed via quantum polylogarithms and their generalizations. 
 In a given cluster coordinate system ${\bf c}$, the asymptotic expansion of the vectors ${\cal A}^\hbar_{n, {\rm L}}$ should have the form  
 \be \la{1/1/2026.b}
{\cal A}^\hbar_{n, {\rm L}} \sim e^{\alpha_{{\bf c}, {\rm L}}(a, \hbar)/\hbar}, \ \ \ \mbox{where 
 $\alpha_{{\bf c}, {\rm L}}(a, \hbar)$ is a function on 
 $
 {\rm Conf}_n(\R_{>0}^4)$}.
 \ee
 One should have 
 $$
  {\alpha}_{n, {\rm L}}(a, \hbar)_{|\hbar =0} = \mbox{\it the  ${\rm L}-$loop scatterring amplitude}.
 $$
 
 Similar  conjectural cluster description of the correlation functions in the Liouville and Toda theories is discussed in \cite[Section 6]{GS}.

\end{document}